\documentclass{amsart}

\usepackage{amssymb,amsmath,amsthm,latexsym,booktabs}

\theoremstyle{plain}

\newtheorem{conjecture}{Conjecture}

\newcommand{\nn}{\!\!}

\begin{document}

\title{Canonical forms of small tensors over $\mathbb{F}_2$}

\author{Murray R. Bremner}

\address{Department of Mathematics and Statistics, University of Saskatchewan, Canada}

\email{bremner@math.usask.ca}

\author{Jiaxiong Hu}

\address{Department of Mathematics, Simon Fraser University, Canada}

\email{hujiaxiong@gmail.com}

\date{20 June 2012}



\begin{abstract}
We consider multidimensional arrays with at most 27 entries over the field with two elements,
and their equivalence classes for the action of the direct product of general linear groups.
The possible 3-dimensional formats are 
$p \times 2 \times 2$ ($p = 2, \dots, 6$), 
$p \times 3 \times 2$ ($p = 3, 4$), 
and
$3 \times 3 \times 3$;
the possible 4-dimensional formats are 
$p \times 2 \times 2 \times 2$ ($p = 2, 3$).
In each case, we compute the orbits for the group action, and then we determine the rank of each orbit.
In particular, we determine the maximum rank for these arrays over $\mathbb{F}_2$.
\end{abstract}

\maketitle

%%%%%%%%%%%%%%%%%%%%%%%%%%%%%%%%%%%%%%%%%%%%%%%%%%%%%%%%%%%%%%%%%%%%%%%%%%%%%%

\section{Introduction}

By a tensor we mean an element of the tensor product of vector spaces.
The rank of a tensor is the minimal number of terms in its expression as a sum of simple tensors.
For a recent survey, with emphasis on algorithms and
applications over fields of characteristic 0, see Kolda and Bader \cite{KB}.
An equivalent problem is the computational complexity of sets of multilinear forms;
see B\"urgisser et al.~\cite{BCS}.

In the classical case of two vector spaces over any field $\mathbb{F}$, we can express 
the problem in terms of $p \times q$ matrices ($p \ge q$).
The orbit representatives for the action of $GL_p(\mathbb{F}) \times GL_q(\mathbb{F})$ 
are the $q+1$ matrices of the form
  \[
  \left[
  \begin{array}{cc}
  I_r & 0 \\
  0 & 0
  \end{array}
  \right]
  \quad
  ( r = 0,\dots, q ),
  \]
where $I_r$ is the $r \times r$ identity matrix and $r$ is the rank.
For more than two vector spaces, the situation is much more complicated.

In this paper, we consider vector spaces over the field $\mathbb{F}_2$ with two elements;
the maximum rank of such a tensor can be higher than expected.
For example, Kruskal \cite{K} observed that over $\mathbb{R}$ the rank of a $3 \times 3 \times 3$
tensor is at most 5, but von zur Gathen \cite{vzG} has shown that over $\mathbb{F}_2$ there exist
$3 \times 3 \times 3$ tensors of rank 6.
In this paper we verify that the maximum rank is 6 and
exhibit three inequivalent canonical forms of this rank.
For tensors of format $p \times q \times 2$, see the recent results by Sumi et al.~\cite{S},
which complete earlier results of Ja'Ja' \cite{J}.
We verify that for tensors of formats $3 \times 3 \times 2$ and $4 \times 3 \times 2$, 
the maximum rank is 5, and we obtain respectively one and four distinct canonical forms.

Limitations on available computer memory impose a maximum of 27 entries, 
which restricts our study to the formats
$p \times 2 \times 2$ ($p = 2, \dots, 6$), 
$p \times 3 \times 2$ ($p = 3, 4$), 
$3 \times 3 \times 3$, and
$p \times 2 \times 2 \times 2$ ($p = 2, 3$).
In each case, the set of tensors admits an action by the direct product of general linear groups.
The canonical form of such a tensor is the lexicographically minimal element of its orbit under this group action.
Since the rank is the same for tensors in the same orbit,
once we have found the canonical forms, it is not difficult to find the ranks.

A summary of our results appears in Table \ref{summary}.
In the first seven cases, the group of symmetries is the ``small group'':
the direct product of general linear groups.
In the last three cases, the group of symmetries is the ``large group'':
the semidirect product of the small group with the permutations
of the equal dimensions.
Orbits for the action of the small (large) group will be called small (large) orbits.

All computations were done with Maple 16 on a
Lenovo ThinkCentre M91p Tower 7052A8U i7-2600 CPU (Quad Core 3.40/3.80GHz)
using Windows 7 Professional 64-bit
with 16 gigabytes of RAM.

\begin{table}
\begin{center}
\begin{tabular}{lrrrr}
format &\qquad $\#$ tensors &\qquad group order & $\#$ orbits & max rank 
\\ \midrule
$2 \times 2 \times 2$ & 256 & 216 & 8 & 3
\\
$3 \times 2 \times 2$ & 4096 & 6048 & 10 & 3
\\
$4 \times 2 \times 2$ & 65536 & 725760 & 11 & 4
\\
$5 \times 2 \times 2$ & 1048576 & 359976960 & 11 & 4 
\\
$6 \times 2 \times 2$ & 16777216 & 725713551360 & 11 & 4
\\
$3 \times 3 \times 2$ & 262144 & 169344 & 21 & 5
\\
$4 \times 3 \times 2$ & 16777216 & 20321280 & 28 & 5
\\ \midrule
$3 \times 3 \times 3$ & 134217728 & 28449792 & 56 & 6
\\
$2 \times 2 \times 2 \times 2$ & 65536 & 31104 & 31 & 6
\\
$3 \times 2 \times 2 \times 2$ & 16777216 & 217728 & 213 & 6
\\ \midrule
\end{tabular}
\end{center}
\caption{Summary of results}
\label{summary}
\end{table}

%%%%%%%%%%%%%%%%%%%%%%%%%%%%%%%%%%%%%%%%%%%%%%%%%%%%%%%%%%%%%%%%%%%%%%%%%%%%%%

\section{Preliminaries}

Let $V_1, \dots, V_n$ be finite dimensional vector spaces over a field $\mathbb{F}$.
By a tensor we mean an element of the tensor product $V_1 \otimes \cdots \otimes V_n$.
A simple tensor is one of the form $v_1 \otimes \cdots \otimes v_n$ where $v_i \in V_i$, $v_i \ne 0$
($i = 1, \dots, n$).
Every tensor is a sum of simple tensors, and the rank of a tensor is the minimal number
of simple tensors occuring in such a representation.
If we choose bases in each $V_i$ then we can identify a tensor $X$ with a multidimensional
array $( x_{i_1 \cdots i_n} )$ of format $d_1 \times \cdots \times d_n$ where 
$d_i = \dim V_i$ ($i = 1, \dots, n$).
We also use the term tensor for such an array, in order to avoid confusion with the data structures
called arrays in Maple.

The flattening of a tensor $X = ( x_{i_1 \cdots i_n} )$ is the row vector
  \[
  \mathrm{flat}(X) = [ x_{1 \cdots 1}, \dots, x_{i_1 \cdots i_n}, \dots, x_{d_1 \cdots d_n} ],
  \]
where the entries are in lex order by subscripts: $i_1 \cdots i_n$ precedes $i'_1 \cdots i'_n$
if and only if $i_j < i'_j$ where $j$ is the least index for which $i_j \ne i'_j$.
Conversely, the unflattening of such a row vector is the corresponding tensor.
If $\mathbb{F}$ is the prime field $\mathbb{F}_m$ with $m$ elements, then we can encode a
tensor $X$ as the non-negative integer whose representation in base $m$ is $\mathrm{flat}(X)$.
Conversely, the decoding of an integer in the range from 0 to $m^{d_1 \cdots d_n}-1$ is the
corresponding tensor.
The lex order on flattenings coincides with the natural order on integers.
The minimal element of a set of tensors is defined in terms of this total order.

The direct product $GL(V_1) \times \cdots \times GL(V_n)$ of general linear groups acts on
$V_1 \otimes \cdots \otimes V_n$ in the natural way.
The canonical form of a tensor is the minimal element in its orbit under this group action.
For $\mathbb{F} = \mathbb{F}_m$, the finite group $GL(V)$, $\dim V = d$, has order
$(m^d-1)(m^d-m) \cdots (m^d-m^{d-1})$.
If we choose a basis $e_1, \dots, e_d$ for $V$ then elements of $GL(V)$ correspond to 
$d \times d$ matrices over $\mathbb{F}_m$.
For $m = 2$, the group of invertible $d \times d$ matrices is generated by two elements: 
the cyclic permutation 
$e_i \mapsto e_{i+1}$ ($i = 1, \dots, d-1$), $e_d \mapsto e_1$, 
and the row operation
$e_1 \mapsto e_1 + e_2$, $e_i \mapsto e_i$ ($i = 2, \dots, d$).
Hence the direct product $GL(V_1) \times \cdots \times GL(V_n)$ is generated by a set 
$\mathcal{G}$ of $2n$ elements.

For a tensor $X$ over $\mathbb{F}_2$, we use the spinning algorithm to compute its orbit.
In the following pseudocode, $\mathcal{O}$ is the current value of the orbit, 
$\mathcal{L}$ contains the new elements computed during the previous iteration,
and $\mathcal{N}$ contains the new elements computed during the current iteration:
  \begin{enumerate}
  \item
  $\mathcal{O} \leftarrow \emptyset$;
  $\mathcal{L} \leftarrow \{ X \}$ 
  \item
  while $\mathcal{L} \ne \emptyset$ do:
    \begin{enumerate}
    \item
    $\mathcal{O} \leftarrow \mathcal{O} \cup \mathcal{L}$      
    \item
    $\mathcal{N} \leftarrow \emptyset$;
    for $Y \in \mathcal{L}$ do for $M \in \mathcal{G}$ do:
      $\mathcal{N} \leftarrow \mathcal{N} \cup \{ M \cdot Y \}$
    \item
    $\mathcal{L} \leftarrow \mathcal{N} \setminus \mathcal{O}$
    \end{enumerate}
  \item
  return $\mathcal{O}$
  \end{enumerate}
We first create a large Maple array, called \texttt{orbitarray}, with $2^N-1$ entries,
where $N = d_1 \cdots d_n$ is the number of entries in the tensors under consideration.  
The indices of \texttt{orbitarray} correspond to nonzero tensors: for an index $i = 1, \dots, 2^N-1$ 
we first decode $i$ by writing it as a binary numeral of $N$ bits (adding leading 0s if necessary), 
and then unflatten this binary numeral to obtain the corresponding tensor.  
To start, every entry of \texttt{orbitarray} is set to 0.
We then perform the following iteration:
  \begin{enumerate}
  \item
  $\omega \leftarrow 0$, $i \leftarrow 0$
  \item
  while $i < 2^N-1$ do:
    \begin{enumerate}
    \item
    $i \leftarrow i + 1$
    \item
    if $\texttt{orbitarray}[i] = 0$ then 
      \begin{enumerate}
      \item
      $\omega \leftarrow \omega + 1$
      \item
      $\texttt{findorbit}[i]$
      \end{enumerate}
    \end{enumerate}
  \end{enumerate}
Procedure \texttt{findorbit} takes the index $i$, decodes and unflattens it 
to the corresponding tensor $X$, 
uses the spinning algorithm to generate the orbit $\mathcal{O}(X)$, and 
sets the corresponding entries of \texttt{orbitarray} to the orbit index $\omega$.
Upon termination, $\omega$ equals the total number of orbits for the group
action, and \texttt{orbitarray} represents the function which assigns to each tensor the index
number of its orbit.
The natural order of the index numbers of the orbits agrees with the lex order 
on the minimal elements in the orbits (the canonical forms of the tensors).

The next step is to compute the ranks of the orbits.
We create another Maple array, called \texttt{linkarray}, of the same size as \texttt{orbitarray}.
We use the data from \texttt{orbitarray} to set entry $i$ of \texttt{linkarray} 
(representing the tensor $X$)
equal to the index $j$ of the next tensor in lex order in the orbit containing $X$.
We then create another Maple array of the same size, called \texttt{rankarray}, 
and initialize every entry to 0.
We generate all simple tensors (tensor products of nonzero vectors) and 
set the corresponding entries of \texttt{rankarray} to 1.
Each index $i$ for which $\texttt{rankarray}[i] = 1$ represents the encoding of a tensor of rank 1.
Let $E$ denote the minimal tensor of rank 1: its flattening is $[0,\dots,0,1]$.
We then perform the following iteration:
  \begin{enumerate}
  \item
  $\texttt{oldrank} \leftarrow 0$, $\texttt{finished} \leftarrow \textrm{false}$
  \item
  While not \texttt{finished} do:
    \begin{enumerate}
    \item
    $\texttt{oldrank} \leftarrow \texttt{oldrank} + 1$,
    $\texttt{finished} \leftarrow \textrm{true}$
    \item
    For each index $i$ for which $\texttt{rankarray}[i] = \texttt{oldrank}$, do:
      \begin{enumerate}
      \item
      Let $X$ be the unflattening of the decoding of $i$.
      \item
      Set $Y \leftarrow X + E$:
      this amounts to changing the rightmost bit of the flattening of $X$ from 0 to 1 or from 1 to 0.
      \item
      Let $j$ be the encoding of the flattening of $Y$.
      Thus $j = i+1$ if $i$ is even, and $j = i-1$ if $i$ is odd.
      \item
      If $\texttt{rankarray}[j] = 0$, then $Y$ has rank $\texttt{oldrank}+1$.
      In this case:
        \begin{itemize}
        \item      
        Use \texttt{linkarray} to store $\texttt{oldrank}+1$ in every entry
        of \texttt{rankarray} corresponding to the tensors in the orbit of $Y$.
        \item
        $\texttt{finished} \leftarrow \textrm{false}$
        \end{itemize}
      \end{enumerate}
    \end{enumerate}
  \end{enumerate}
The iteration terminates when every entry of \texttt{rankarray} contains a positive integer,
which is the rank of the corresponding (nonzero) tensor.

%%%%%%%%%%%%%%%%%%%%%%%%%%%%%%%%%%%%%%%%%%%%%%%%%%%%%%%%%%%%%%%%%%%%%%%%%%%%%%

\section{Formats $p \times 2 \times 2$ ($p = 2, \dots, 6$)}

For $p = 2$ there are $2^8 = 256$ tensors, and the group order is $6^3 = 216$.
There are 7 nonzero orbits with maximum rank 3.  Here is a summary:
\[
\begin{array}{lrrrr}
\text{rank} 
&\quad 0 &\quad 1 &\quad 2 &\quad 3
\\
\text{$\#$ orbits} 
&\quad 1 &\quad 1 &\quad 4 &\quad 2
\\
\text{$\#$ tensors} 
&\quad 1 &\quad 27 &\quad 162 &\quad 66
\\
\text{percent}
&\quad 0.3906
&\quad 10.5469
&\quad 63.2813
&\quad 25.7813
\end{array}
\]
For the orbit sizes and canonical forms, see Table \ref{table222};
the orbits are sorted first by increasing rank, within each rank by increasing orbit size, 
and within each orbit size by lex order of the canonical forms.
This computation took less than 0.1 second. 
A similar classification appears in \cite{BS}; however, that work includes the symmetries
obtained by permuting the three directions.

For $p = 3$ there are $2^{12} = 4096$ tensors, and the group order is $168 \cdot 6^2 = 6048$.
There are 9 nonzero orbits with maximum rank 3:
\[
\begin{array}{lrrrr}
\text{rank} 
&\quad 0 &\quad 1 &\quad 2 &\quad 3
\\
\text{$\#$ orbits} 
&\quad 1 &\quad 1 &\quad 4 &\quad 4
\\
\text{$\#$ tensors} 
&\quad 1 &\quad 63 &\quad 1050 &\quad 2982
\\
\text{percent}
&\quad 0.0244
&\quad 1.5381
&\quad 25.6348
&\quad 72.8027
\end{array}
\]
For the orbit sizes and canonical forms, see Table \ref{table322}.
This computation took less than 1 second. 

For $p = 4$ there are $2^{16} = 65536$ tensors, and the group order is $20160 \cdot 6^2 = 725760$.
There are 10 nonzero orbits with maximum rank 4:
\[
\begin{array}{lrrrrr}
\text{rank} 
&\quad 0 &\quad 1 &\quad 2 &\quad 3 &\quad 4
\\
\text{$\#$ orbits} 
&\quad 1 &\quad 1 &\quad 4 &\quad 4 &\quad 1
\\
\text{$\#$ tensors} 
&\quad 1 &\quad 135 &\quad 5130 &\quad 40110 &\quad 20160
\\
\text{percent}
&\quad 0.0015
&\quad 0.2060
&\quad 7.8278
&\quad 61.2030
&\quad 30.7617
\end{array}
\]
For the orbit sizes and canonical forms, see Table \ref{table422}.
One of the orbits of rank 3 contains almost 35\% of the tensors; 
the single orbit of rank 4 contains almost 31\% of the tensors.
This computation took less than 10 seconds. 

\begin{table}
\begin{center}
\begin{tabular}{rrrcccccccc}
$\#$ & rank & size & \multicolumn{8}{c}{canonical form} \\ \midrule
    1 & 1 &      27 & .& .& .& .& .& .& .& 1 \\ \midrule 
    2 & 2 &      18 & .& .& .& .& .& 1& 1& . \\ 
    3 & 2 &      18 & .& .& .& 1& .& .& 1& . \\ 
    4 & 2 &      18 & .& .& .& 1& .& 1& .& . \\ 
    5 & 2 &     108 & .& .& .& 1& 1& .& .& . \\ \midrule 
    6 & 3 &      12 & .& 1& 1& .& 1& .& 1& 1 \\ 
    7 & 3 &      54 & .& .& .& 1& .& 1& 1& . \\ \midrule
\end{tabular}
\end{center}
\caption{Small orbits of $2 \times 2 \times 2$ tensors}
\label{table222}
\bigskip
\begin{center}
\begin{tabular}{rrrcccccccccccc}
$\#$ & rank & size & \multicolumn{12}{c}{canonical form} \\ \midrule
    1 & 1 &      63 & .& .& .& .& .& .& .& .& .& .& .& 1 \\ \midrule
    2 & 2 &      42 & .& .& .& .& .& .& .& .& .& 1& 1& . \\ 
    3 & 2 &     126 & .& .& .& .& .& .& .& 1& .& .& 1& . \\ 
    4 & 2 &     126 & .& .& .& .& .& .& .& 1& .& 1& .& . \\ 
    5 & 2 &     756 & .& .& .& .& .& .& .& 1& 1& .& .& . \\ \midrule
    6 & 3 &      84 & .& .& .& .& .& 1& 1& .& 1& .& 1& 1 \\ 
    7 & 3 &     378 & .& .& .& .& .& .& .& 1& .& 1& 1& . \\ 
    8 & 3 &    1008 & .& .& .& 1& .& 1& 1& .& 1& .& .& . \\ 
    9 & 3 &    1512 & .& .& .& 1& .& .& 1& .& .& 1& .& . \\ \midrule
\end{tabular}
\end{center}
\caption{Small orbits of $3 \times 2 \times 2$ tensors}
\label{table322}
\bigskip
\begin{center}
\begin{tabular}{rrrcccccccccccccccc}
$\#$ & rank & size & \multicolumn{16}{c}{canonical form} \\ \midrule
    1 & 1 &     135 & .& .& .& .& .& .& .& .& .& .& .& .& .& .& .& 1 \\ \midrule 
    2 & 2 &      90 & .& .& .& .& .& .& .& .& .& .& .& .& .& 1& 1& . \\ 
    3 & 2 &     630 & .& .& .& .& .& .& .& .& .& .& .& 1& .& .& 1& . \\ 
    4 & 2 &     630 & .& .& .& .& .& .& .& .& .& .& .& 1& .& 1& .& . \\ 
    5 & 2 &    3780 & .& .& .& .& .& .& .& .& .& .& .& 1& 1& .& .& . \\ \midrule 
    6 & 3 &     420 & .& .& .& .& .& .& .& .& .& 1& 1& .& 1& .& 1& 1 \\ 
    7 & 3 &    1890 & .& .& .& .& .& .& .& .& .& .& .& 1& .& 1& 1& . \\ 
    8 & 3 &   15120 & .& .& .& .& .& .& .& 1& .& 1& 1& .& 1& .& .& . \\ 
    9 & 3 &   22680 & .& .& .& .& .& .& .& 1& .& .& 1& .& .& 1& .& . \\ \midrule 
   10 & 4 &   20160 & .& .& .& 1& .& .& 1& .& .& 1& .& .& 1& .& .& . \\ \midrule
\end{tabular}
\end{center}
\caption{Small orbits of $4 \times 2 \times 2$ tensors}
\label{table422}
\end{table}

For $p = 5$ and $p = 6$ we obtain the same canonical forms with the addition of more leading 0s
in the flattenings of the tensors; thus in both cases there are 10 nonzero orbits with maximum rank 4.

For $p = 5$ there are $2^{20} = 1048576$ tensors; the group order is $9999360 \cdot 6^2 = 359976960$.
Here is a summary:
\[
\begin{array}{lrrrrr}
\text{rank} 
&\quad 0 &\quad 1 &\quad 2 &\quad 3 &\quad 4
\\
\text{$\#$ orbits} 
&\quad 1 &\quad 1 &\quad 4 &\quad 4 &\quad 1
\\
\text{$\#$ tensors} 
&\quad 1 &\quad 279 &\quad 22506 &\quad 400830 &\quad 624960
\\
\text{percent}
&\quad 0.0001
&\quad 0.0266
&\quad 2.1463
&\quad 38.2261
&\quad 59.6008
\end{array}
\]
The orbit sizes, using the order of Table \ref{table422} (with four more leading 0s), are:
  \[ 
     279, \; 
     186, \;  
    2790, \;  
    2790, \; 
   16740, \;  
    1860, \;  
    8370, \;  
  156240, \;  
  234360, \;  
  624960.
  \]  
The single orbit of rank 4 contains almost 60\% of the tensors.
This computation took less than 4 minutes. 

For $p = 6$ there are $2^{24} = 16777216$ tensors; 
the group order is $20158709760 \cdot 6^2 = 725713551360$.
Here is a summary:
\[
\begin{array}{lrrrrr}
\text{rank} 
&\quad 0 &\quad 1 &\quad 2 &\quad 3 &\quad 4
\\
\text{$\#$ orbits} 
&\quad 1 &\quad 1 &\quad 4 &\quad 4 &\quad 1
\\
\text{$\#$ tensors} 
&\quad 1 &\quad 567 &\quad 94122 &\quad 3558366 &\quad 13124160
\\
\text{percent}
&\quad 0.0000
&\quad 0.0034
&\quad 0.5610
&\quad 21.2095
&\quad 78.2261
\end{array}
\]
The orbit sizes, using the order of Table \ref{table422} (with eight more leading 0s), are:
  \[ 
       567, \; 
       378, \; 
     11718, \; 
     11718, \; 
     70308, \; 
      7812, \; 
     35154, \; 
   1406160, \; 
   2109240, \; 
  13124160. 
  \]  
The single orbit of rank 4 contains more than 78\% of the tensors.
This computation took just over 75 minutes. 

These results suggest the following conjectures for $p \times 2 \times 2$ tensors with $p \ge 4$.

\begin{conjecture} \label{conjecturep22}
For $p \ge 4$, there are 10 canonical forms for nonzero $p \times 2 \times 2$ tensors over $\mathbb{F}_2$.
For each orbit, the rank and the horizontal $2 \times 2$ slices of the minimal element
(obtained by fixing the first subscript) are given in Table \ref{tablep22}.
The leftmost empty slice indicates the appropriate number of zero $2 \times 2$ slices.
\end{conjecture}

\begin{table}
\begin{center}
\begin{tabular}{rrr}
$\#$ & rank &\quad canonical form \\ \midrule
  1 & 1 &\quad
  $\left[
  \begin{array}{c|cc|cc|cc|cc}
  \;\;\;\;\; & \cdot & \cdot & \cdot & \cdot & \cdot & \cdot & \cdot & \cdot \\
  \;\;\;\;\; & \cdot & \cdot & \cdot & \cdot & \cdot & \cdot & \cdot & 1 
  \end{array}
  \right]$
  \\ \midrule 
  2 & 2 &\quad 
  $\left[
  \begin{array}{c|cc|cc|cc|cc}
  \;\;\;\;\; & \cdot & \cdot & \cdot & \cdot & \cdot & \cdot & \cdot & 1 \\
  \;\;\;\;\; & \cdot & \cdot & \cdot & \cdot & \cdot & \cdot & 1 & \cdot 
  \end{array}
  \right]$
  \\[8pt]    
  3 & 2 &\quad 
  $\left[
  \begin{array}{c|cc|cc|cc|cc}
  \;\;\;\;\; & \cdot & \cdot & \cdot & \cdot & \cdot & \cdot & \cdot & \cdot \\
  \;\;\;\;\; & \cdot & \cdot & \cdot & \cdot & \cdot & 1 & 1 & \cdot
  \end{array}
  \right]$
  \\[8pt]   
  4 & 2 &\quad 
  $\left[
  \begin{array}{c|cc|cc|cc|cc}
  \;\;\;\;\; & \cdot & \cdot & \cdot & \cdot & \cdot & \cdot & \cdot & 1 \\
  \;\;\;\;\; & \cdot & \cdot & \cdot & \cdot & \cdot & 1 & \cdot & \cdot 
  \end{array}
  \right]$
  \\[8pt]    
  5 & 2 &\quad 
  $\left[
  \begin{array}{c|cc|cc|cc|cc}
  \;\;\;\;\; & \cdot & \cdot & \cdot & \cdot & \cdot & \cdot & 1 & \cdot \\
  \;\;\;\;\; & \cdot & \cdot & \cdot & \cdot & \cdot & 1 & \cdot & \cdot 
  \end{array}
  \right]$
  \\ \midrule 
  6 & 3 &\quad 
  $\left[
  \begin{array}{c|cc|cc|cc|cc}
  \;\;\;\;\; & \cdot & \cdot & \cdot & \cdot & \cdot & 1 & 1 & \cdot \\
  \;\;\;\;\; & \cdot & \cdot & \cdot & \cdot & 1 & \cdot & 1 & 1 
  \end{array}
  \right]$
  \\[8pt]    
  7 & 3 &\quad 
  $\left[
  \begin{array}{c|cc|cc|cc|cc}
  \;\;\;\;\; & \cdot & \cdot & \cdot & \cdot & \cdot & \cdot & \cdot & 1 \\
  \;\;\;\;\; & \cdot & \cdot & \cdot & \cdot & \cdot & 1 & 1 & \cdot 
  \end{array}
  \right]$
  \\[8pt]    
  8 & 3 &\quad 
  $\left[
  \begin{array}{c|cc|cc|cc|cc}
  \;\;\;\;\; & \cdot & \cdot & \cdot & \cdot & \cdot & 1 & 1 & \cdot \\
  \;\;\;\;\; & \cdot & \cdot & \cdot & 1 & 1 & \cdot & \cdot & \cdot 
  \end{array}
  \right]$
  \\[8pt]    
  9 & 3 &\quad 
  $\left[
  \begin{array}{c|cc|cc|cc|cc}
  \;\;\;\;\; & \cdot & \cdot & \cdot & \cdot & \cdot & \cdot & \cdot & 1 \\
  \;\;\;\;\; & \cdot & \cdot & \cdot & 1 & 1 & \cdot & \cdot & \cdot 
  \end{array}
  \right]$
  \\ \midrule      
  10 & 4 &\quad 
  $\left[
  \begin{array}{c|cc|cc|cc|cc}
  \;\;\;\;\; & \cdot & \cdot & \cdot & \cdot & \cdot & 1 & 1 & \cdot \\
  \;\;\;\;\; & \cdot & 1 & 1 & \cdot & \cdot & \cdot & \cdot & \cdot 
  \end{array}
  \right]$
  \\ \midrule    
\end{tabular}
\end{center}
\caption{Canonical forms for Conjecture \ref{conjecturep22}}
\label{tablep22}
\end{table}

\begin{conjecture} \label{conjecturep22x}
As $p$ becomes arbitrarily large, the set of $p \times 2 \times 2$ tensors over $\mathbb{F}_2$
contains a single orbit $\mathcal{O}$ of rank 4 which contains almost all of the $2^{4p}$ tensors:
  \[
  \lim_{p \to \infty} \frac{|\mathcal{O}|}{2^{4p}} = 1.
  \]
\end{conjecture}

%%%%%%%%%%%%%%%%%%%%%%%%%%%%%%%%%%%%%%%%%%%%%%%%%%%%%%%%%%%%%%%%%%%%%%%%%%%%%%

\section{Formats $p \times 3 \times 2$ ($p = 3, 4$)}

For $p = 3$ there are $2^{18} = 262144$ tensors, and the group order is $168^2 \cdot 6 = 169344$.
There are 20 nonzero orbits with maximum rank 5:
\[
\begin{array}{lrrrrrr}
\text{rank} 
&\quad 0 &\quad 1 &\quad 2 &\quad 3 &\quad 4 &\quad 5
\\
\text{$\#$ orbits} 
&\quad 1 &\quad 1 &\quad 4 &\quad 9 &\quad 5 &\quad 1
\\
\text{$\#$ tensors} 
&\quad 1 &\quad 147 &\quad 6762 &\quad 95466 &\quad 151704 &\quad 8064
\\
\text{percent}
&\quad  0.0004
&\quad  0.0561
&\quad  2.5795
&\quad 36.4174
&\quad 57.8705
&\quad  3.0762
\end{array}
\]
For the orbit sizes and canonical forms, see Table \ref{table332}.
This computation took less than 30 seconds. 

\begin{table}
\begin{center}
\begin{tabular}{rrrcccccccccccccccccc}
$\#$ & rank & size & \multicolumn{18}{c}{canonical form} \\ \midrule
   1 & 1 &   147 & .& .& .& .& .& .& .& .& .& .& .& .& .& .& .& .& .& 1 \\ \midrule 
   2 & 2 &   294 & .& .& .& .& .& .& .& .& .& .& .& .& .& .& .& 1& 1& . \\ 
   3 & 2 &   294 & .& .& .& .& .& .& .& .& .& .& .& 1& .& .& .& .& 1& . \\ 
   4 & 2 &   882 & .& .& .& .& .& .& .& .& .& .& .& 1& .& .& .& 1& .& . \\ 
   5 & 2 &  5292 & .& .& .& .& .& .& .& .& .& .& .& 1& .& .& 1& .& .& . \\ \midrule 
   6 & 3 &   504 & .& .& .& .& .& 1& .& .& .& 1& .& .& .& 1& .& .& .& . \\ 
   7 & 3 &   588 & .& .& .& .& .& .& .& .& .& 1& 1& .& .& .& 1& .& 1& 1 \\ 
   8 & 3 &  2646 & .& .& .& .& .& .& .& .& .& .& .& 1& .& .& .& 1& 1& . \\ 
   9 & 3 &  7056 & .& .& .& .& .& .& .& .& .& 1& 1& .& .& 1& 1& .& .& . \\ 
  10 & 3 &  7056 & .& .& .& .& .& 1& .& .& .& 1& 1& .& .& .& 1& .& .& . \\ 
  11 & 3 & 10584 & .& .& .& .& .& .& .& .& .& .& .& 1& .& 1& 1& .& .& . \\ 
  12 & 3 & 10584 & .& .& .& .& .& 1& .& .& .& .& 1& .& .& .& .& 1& .& . \\ 
  13 & 3 & 28224 & .& .& .& .& .& 1& .& .& .& 1& .& .& 1& .& .& .& .& . \\ 
  14 & 3 & 28224 & .& .& .& .& .& 1& .& .& 1& .& .& .& 1& 1& .& .& .& . \\ \midrule 
  15 & 4 &  7056 & .& .& .& .& .& 1& .& .& .& .& 1& .& .& 1& 1& .& .& . \\ 
  16 & 4 & 10584 & .& .& .& .& .& 1& .& .& .& 1& .& .& .& 1& .& .& 1& . \\ 
  17 & 4 & 21168 & .& .& .& .& .& 1& .& .& .& 1& 1& .& .& 1& 1& .& .& . \\ 
  18 & 4 & 28224 & .& .& .& .& .& 1& .& 1& 1& .& .& .& 1& .& 1& 1& .& . \\ 
  19 & 4 & 84672 & .& .& .& .& .& 1& .& .& .& 1& 1& .& 1& .& .& .& .& . \\ \midrule 
  20 & 5 &  8064 & .& .& .& 1& 1& .& .& 1& 1& .& .& .& 1& .& .& .& 1& 1 \\ \midrule 
\end{tabular}
\end{center}
\caption{Small orbits of $3 \times 3 \times 2$ tensors}
\label{table332}
\end{table}

For $p = 4$ there are $2^{24} = 16777216$ tensors, and the group order is 
$20160 \cdot 168 \cdot 6 = 20321280$.
There are 27 nonzero orbits with maximum rank 5:
\[
\begin{array}{lrrrrrr}
\text{rank} 
&\;\; 0 &\;\; 1 &\;\; 2 &\;\; 3 &\;\; 4 &\;\; 5
\\
\text{$\#$ orbits} 
&\;\; 1 &\;\; 1 &\;\; 4 &\;\; 9 &\;\; 9 &\;\; 4
\\
\text{$\#$ tensors} 
&\;\; 1 &\;\; 315 &\;\; 32970 &\;\; 1223250 &\;\; 10460520 &\;\; 5060160
\\
\text{percent}
&\;\;  0.0000
&\;\;  0.0019
&\;\;  0.1965
&\;\;  7.2911
&\;\; 62.3496
&\;\; 30.1609
\end{array}
\]
For the orbit sizes and canonical forms, see Table \ref{table432}.
This computation took less than 45 minutes. 

\begin{table}
\begin{center}
\begin{tabular}{rrrcccccccccccccccccccccccc}
$\#$ & rank & size & \multicolumn{24}{c}{canonical form} \\ \midrule
    1 & 1 &     315 & .&\nn .&\nn .&\nn .&\nn .&\nn .&\nn .&\nn .&\nn .&\nn .&\nn .&\nn .&\nn .&\nn .&\nn .&\nn .&\nn .&\nn .&\nn .&\nn .&\nn .&\nn .&\nn .&\nn 1 \\ \midrule 
    2 & 2 &     630 & .&\nn .&\nn .&\nn .&\nn .&\nn .&\nn .&\nn .&\nn .&\nn .&\nn .&\nn .&\nn .&\nn .&\nn .&\nn .&\nn .&\nn .&\nn .&\nn .&\nn .&\nn 1&\nn 1&\nn . \\ 
    3 & 2 &    1470 & .&\nn .&\nn .&\nn .&\nn .&\nn .&\nn .&\nn .&\nn .&\nn .&\nn .&\nn .&\nn .&\nn .&\nn .&\nn .&\nn .&\nn 1&\nn .&\nn .&\nn .&\nn .&\nn 1&\nn . \\ 
    4 & 2 &    4410 & .&\nn .&\nn .&\nn .&\nn .&\nn .&\nn .&\nn .&\nn .&\nn .&\nn .&\nn .&\nn .&\nn .&\nn .&\nn .&\nn .&\nn 1&\nn .&\nn .&\nn .&\nn 1&\nn .&\nn . \\ 
    5 & 2 &   26460 & .&\nn .&\nn .&\nn .&\nn .&\nn .&\nn .&\nn .&\nn .&\nn .&\nn .&\nn .&\nn .&\nn .&\nn .&\nn .&\nn .&\nn 1&\nn .&\nn .&\nn 1&\nn .&\nn .&\nn . \\ \midrule 
    6 & 3 &    2940 & .&\nn .&\nn .&\nn .&\nn .&\nn .&\nn .&\nn .&\nn .&\nn .&\nn .&\nn .&\nn .&\nn .&\nn .&\nn 1&\nn 1&\nn .&\nn .&\nn .&\nn 1&\nn .&\nn 1&\nn 1 \\ 
    7 & 3 &    7560 & .&\nn .&\nn .&\nn .&\nn .&\nn .&\nn .&\nn .&\nn .&\nn .&\nn .&\nn 1&\nn .&\nn .&\nn .&\nn 1&\nn .&\nn .&\nn .&\nn 1&\nn .&\nn .&\nn .&\nn . \\ 
    8 & 3 &   13230 & .&\nn .&\nn .&\nn .&\nn .&\nn .&\nn .&\nn .&\nn .&\nn .&\nn .&\nn .&\nn .&\nn .&\nn .&\nn .&\nn .&\nn 1&\nn .&\nn .&\nn .&\nn 1&\nn 1&\nn . \\ 
    9 & 3 &   35280 & .&\nn .&\nn .&\nn .&\nn .&\nn .&\nn .&\nn .&\nn .&\nn .&\nn .&\nn .&\nn .&\nn .&\nn .&\nn 1&\nn 1&\nn .&\nn .&\nn 1&\nn 1&\nn .&\nn .&\nn . \\ 
   10 & 3 &   52920 & .&\nn .&\nn .&\nn .&\nn .&\nn .&\nn .&\nn .&\nn .&\nn .&\nn .&\nn .&\nn .&\nn .&\nn .&\nn .&\nn .&\nn 1&\nn .&\nn 1&\nn 1&\nn .&\nn .&\nn . \\ 
   11 & 3 &  105840 & .&\nn .&\nn .&\nn .&\nn .&\nn .&\nn .&\nn .&\nn .&\nn .&\nn .&\nn 1&\nn .&\nn .&\nn .&\nn 1&\nn 1&\nn .&\nn .&\nn .&\nn 1&\nn .&\nn .&\nn . \\ 
   12 & 3 &  158760 & .&\nn .&\nn .&\nn .&\nn .&\nn .&\nn .&\nn .&\nn .&\nn .&\nn .&\nn 1&\nn .&\nn .&\nn .&\nn .&\nn 1&\nn .&\nn .&\nn .&\nn .&\nn 1&\nn .&\nn . \\ 
   13 & 3 &  423360 & .&\nn .&\nn .&\nn .&\nn .&\nn .&\nn .&\nn .&\nn .&\nn .&\nn .&\nn 1&\nn .&\nn .&\nn .&\nn 1&\nn .&\nn .&\nn 1&\nn .&\nn .&\nn .&\nn .&\nn . \\ 
   14 & 3 &  423360 & .&\nn .&\nn .&\nn .&\nn .&\nn .&\nn .&\nn .&\nn .&\nn .&\nn .&\nn 1&\nn .&\nn .&\nn 1&\nn .&\nn .&\nn .&\nn 1&\nn 1&\nn .&\nn .&\nn .&\nn . \\ \midrule 
   15 & 4 &  105840 & .&\nn .&\nn .&\nn .&\nn .&\nn .&\nn .&\nn .&\nn .&\nn .&\nn .&\nn 1&\nn .&\nn .&\nn .&\nn .&\nn 1&\nn .&\nn .&\nn 1&\nn 1&\nn .&\nn .&\nn . \\ 
   16 & 4 &  141120 & .&\nn .&\nn .&\nn .&\nn .&\nn 1&\nn .&\nn .&\nn .&\nn .&\nn 1&\nn .&\nn .&\nn .&\nn .&\nn 1&\nn .&\nn .&\nn .&\nn .&\nn 1&\nn .&\nn .&\nn . \\ 
   17 & 4 &  158760 & .&\nn .&\nn .&\nn .&\nn .&\nn .&\nn .&\nn .&\nn .&\nn .&\nn .&\nn 1&\nn .&\nn .&\nn .&\nn 1&\nn .&\nn .&\nn .&\nn 1&\nn .&\nn .&\nn 1&\nn . \\ 
   18 & 4 &  317520 & .&\nn .&\nn .&\nn .&\nn .&\nn .&\nn .&\nn .&\nn .&\nn .&\nn .&\nn 1&\nn .&\nn .&\nn .&\nn 1&\nn 1&\nn .&\nn .&\nn 1&\nn 1&\nn .&\nn .&\nn . \\ 
   19 & 4 &  423360 & .&\nn .&\nn .&\nn .&\nn .&\nn .&\nn .&\nn .&\nn .&\nn .&\nn .&\nn 1&\nn .&\nn 1&\nn 1&\nn .&\nn .&\nn .&\nn 1&\nn .&\nn 1&\nn 1&\nn .&\nn . \\ 
   20 & 4 &  423360 & .&\nn .&\nn .&\nn .&\nn .&\nn 1&\nn .&\nn .&\nn .&\nn .&\nn 1&\nn .&\nn .&\nn .&\nn .&\nn 1&\nn .&\nn .&\nn .&\nn 1&\nn .&\nn .&\nn .&\nn . \\ 
   21 & 4 & 1270080 & .&\nn .&\nn .&\nn .&\nn .&\nn .&\nn .&\nn .&\nn .&\nn .&\nn .&\nn 1&\nn .&\nn .&\nn .&\nn 1&\nn 1&\nn .&\nn 1&\nn .&\nn .&\nn .&\nn .&\nn . \\ 
   22 & 4 & 2540160 & .&\nn .&\nn .&\nn .&\nn .&\nn 1&\nn .&\nn .&\nn .&\nn .&\nn 1&\nn .&\nn .&\nn .&\nn .&\nn 1&\nn .&\nn .&\nn 1&\nn .&\nn .&\nn .&\nn .&\nn . \\ 
   23 & 4 & 5080320 & .&\nn .&\nn .&\nn .&\nn .&\nn 1&\nn .&\nn .&\nn .&\nn 1&\nn .&\nn .&\nn .&\nn 1&\nn .&\nn .&\nn 1&\nn .&\nn 1&\nn .&\nn .&\nn .&\nn .&\nn . \\ \midrule 
   24 & 5 &  120960 & .&\nn .&\nn .&\nn .&\nn .&\nn .&\nn .&\nn .&\nn .&\nn 1&\nn 1&\nn .&\nn .&\nn 1&\nn 1&\nn .&\nn .&\nn .&\nn 1&\nn .&\nn .&\nn .&\nn 1&\nn 1 \\ 
   25 & 5 &  282240 & .&\nn .&\nn .&\nn .&\nn .&\nn 1&\nn .&\nn .&\nn .&\nn .&\nn 1&\nn .&\nn .&\nn 1&\nn 1&\nn .&\nn .&\nn .&\nn 1&\nn .&\nn 1&\nn 1&\nn .&\nn . \\ 
   26 & 5 & 1270080 & .&\nn .&\nn .&\nn .&\nn .&\nn 1&\nn .&\nn .&\nn .&\nn .&\nn 1&\nn .&\nn .&\nn .&\nn .&\nn 1&\nn .&\nn .&\nn .&\nn 1&\nn 1&\nn .&\nn .&\nn . \\ 
   27 & 5 & 3386880 & .&\nn .&\nn .&\nn .&\nn .&\nn 1&\nn .&\nn .&\nn .&\nn 1&\nn 1&\nn .&\nn .&\nn 1&\nn 1&\nn .&\nn .&\nn .&\nn 1&\nn .&\nn .&\nn .&\nn .&\nn . \\ \midrule 
\end{tabular}
\end{center}
\caption{Small orbits of $4 \times 3 \times 2$ tensors}
\label{table432}
\end{table}

%%%%%%%%%%%%%%%%%%%%%%%%%%%%%%%%%%%%%%%%%%%%%%%%%%%%%%%%%%%%%%%%%%%%%%%%%%%%%%

\section{Format $3 \times 3 \times 3$}

In this and the following section, the formats have three (or four) equal dimensions.
We exploit this symmetry in order to reduce the number of orbits
to a manageable size.
For $3 \times 3 \times 3$ tensors, the large group acting on the tensors is
  \[
  G = \big( GL(3,\mathbb{F}_2) \times GL(3,\mathbb{F}_2) \times GL(3,\mathbb{F}_2) \big) \rtimes S_3,
  \]
where the symmetric group $S_3$ permutes the three directions.
As in the previous sections, we first compute the small orbits obtained by the action of 
the direct product of general linear groups.
We then apply the permutations to determine which small orbits combine to make a
single large orbit under the action of $G$.
Given the canonical form $X$ of a small orbit $\mathcal{O}$ with index number $i$, 
we apply the elements of $S_3$ to obtain tensors $X_1 = X, \dots, X_6$.
We then use the Maple arrays, which we have already computed, to find the index numbers
$i_1 = i, \dots, i_6$ of the small orbits containing these tensors.
We conclude that the union $\mathcal{O}_{i_1} \cup \cdots \cup \mathcal{O}_{i_6}$
is a large orbit for the action of $G$.
The canonical form of the tensors in this large orbit is the smallest (in lex order)
of the canonical forms of $\mathcal{O}_{i_1}, \dots, \mathcal{O}_{i_6}$.

With format $3 \times 3 \times 3$ there are $2^{27} = 134217728$ tensors.
The small group order is $168^3 = 4741632$,
and the large group order is $6 \cdot 168^3 = 28449792$.
There are 115 (nonzero) small orbits and 55 (nonzero) large orbits, with maximum rank 6:
\[
\begin{array}{lrrrrrrr}
\text{rank} 
& 0 & 1 & 2 & 3 & 4 & 5 & 6 
\\
\text{$\#$ small} 
& 1 & 1 & 4 & 18 & 44 & 45 & 3
\\
\text{$\#$ large} 
& 1 & 1 & 2 & 8 & 18 & 23 & 3
\\
\text{$\#$ tensors} 
& 1 & 343 & 43218 & 2372286 & 47506872 & 83670048 & 624960 
\\
\text{percent}
& 0.0000
& 0.0003
& 0.0322
& 1.7675
& 35.3954
& 62.3390
& 0.4656
\end{array}
\]
For the large orbit sizes and canonical forms, see Table \ref{table333}.
This computation took just under 282 minutes.

\renewcommand{\nn}{\!\!\!\!}

\newcommand{\sss}{-2pt}

\begin{table} \small
\begin{tabular}{rrrccccccccccccccccccccccccccc}
$\#$ & rank & size & \multicolumn{27}{c}{canonical form} \\ \midrule
   1 & 1 &      343 & .&\nn .&\nn .&\nn .&\nn .&\nn .&\nn .&\nn .&\nn .&\nn .&\nn .&\nn .&\nn .&\nn .&\nn .&\nn .&\nn .&\nn .&\nn .&\nn .&\nn .&\nn .&\nn .&\nn .&\nn .&\nn .&\nn 1 \\ \midrule 
   2 & 2 &     6174 & .&\nn .&\nn .&\nn .&\nn .&\nn .&\nn .&\nn .&\nn .&\nn .&\nn .&\nn .&\nn .&\nn .&\nn .&\nn .&\nn .&\nn .&\nn .&\nn .&\nn .&\nn .&\nn .&\nn 1&\nn .&\nn 1&\nn . \\[\sss]
   3 & 2 &    37044 & .&\nn .&\nn .&\nn .&\nn .&\nn .&\nn .&\nn .&\nn .&\nn .&\nn .&\nn .&\nn .&\nn .&\nn .&\nn .&\nn .&\nn 1&\nn .&\nn .&\nn .&\nn .&\nn 1&\nn .&\nn .&\nn .&\nn . \\ \midrule 
   4 & 3 &     3528 & .&\nn .&\nn .&\nn .&\nn .&\nn .&\nn .&\nn .&\nn .&\nn .&\nn .&\nn .&\nn .&\nn .&\nn .&\nn .&\nn .&\nn .&\nn .&\nn .&\nn 1&\nn .&\nn 1&\nn .&\nn 1&\nn .&\nn . \\[\sss] 
   5 & 3 &     4116 & .&\nn .&\nn .&\nn .&\nn .&\nn .&\nn .&\nn .&\nn .&\nn .&\nn .&\nn .&\nn .&\nn .&\nn 1&\nn .&\nn 1&\nn .&\nn .&\nn .&\nn .&\nn .&\nn 1&\nn .&\nn .&\nn 1&\nn 1 \\[\sss] 
   6 & 3 &    18522 & .&\nn .&\nn .&\nn .&\nn .&\nn .&\nn .&\nn .&\nn .&\nn .&\nn .&\nn .&\nn .&\nn .&\nn .&\nn .&\nn .&\nn 1&\nn .&\nn .&\nn .&\nn .&\nn .&\nn 1&\nn .&\nn 1&\nn . \\[\sss] 
   7 & 3 &   148176 & .&\nn .&\nn .&\nn .&\nn .&\nn .&\nn .&\nn .&\nn .&\nn .&\nn .&\nn .&\nn .&\nn .&\nn 1&\nn .&\nn 1&\nn .&\nn .&\nn .&\nn .&\nn .&\nn 1&\nn .&\nn 1&\nn .&\nn . \\[\sss] 
   8 & 3 &   222264 & .&\nn .&\nn .&\nn .&\nn .&\nn .&\nn .&\nn .&\nn .&\nn .&\nn .&\nn .&\nn .&\nn .&\nn .&\nn .&\nn .&\nn 1&\nn .&\nn .&\nn .&\nn .&\nn 1&\nn .&\nn 1&\nn .&\nn . \\[\sss] 
   9 & 3 &   592704 & .&\nn .&\nn .&\nn .&\nn .&\nn .&\nn .&\nn .&\nn .&\nn .&\nn .&\nn .&\nn .&\nn .&\nn .&\nn .&\nn .&\nn 1&\nn .&\nn 1&\nn .&\nn 1&\nn .&\nn .&\nn .&\nn .&\nn . \\[\sss] 
  10 & 3 &   592704 & .&\nn .&\nn .&\nn .&\nn .&\nn .&\nn .&\nn .&\nn .&\nn .&\nn .&\nn .&\nn .&\nn .&\nn 1&\nn .&\nn 1&\nn .&\nn 1&\nn .&\nn .&\nn .&\nn .&\nn .&\nn .&\nn 1&\nn . \\[\sss] 
  11 & 3 &   790272 & .&\nn .&\nn .&\nn .&\nn .&\nn .&\nn .&\nn .&\nn 1&\nn .&\nn .&\nn .&\nn .&\nn 1&\nn .&\nn .&\nn .&\nn .&\nn 1&\nn .&\nn .&\nn .&\nn .&\nn .&\nn .&\nn .&\nn . \\ \midrule 
  12 & 4 &   148176 & .&\nn .&\nn .&\nn .&\nn .&\nn .&\nn .&\nn .&\nn .&\nn .&\nn .&\nn .&\nn .&\nn .&\nn 1&\nn .&\nn 1&\nn .&\nn .&\nn .&\nn 1&\nn .&\nn .&\nn .&\nn 1&\nn .&\nn . \\[\sss] 
  13 & 4 &   197568 & .&\nn .&\nn .&\nn .&\nn .&\nn 1&\nn .&\nn 1&\nn .&\nn .&\nn .&\nn 1&\nn .&\nn .&\nn .&\nn 1&\nn .&\nn .&\nn .&\nn 1&\nn .&\nn 1&\nn 1&\nn .&\nn .&\nn 1&\nn . \\[\sss] 
  14 & 4 &   222264 & .&\nn .&\nn .&\nn .&\nn .&\nn .&\nn .&\nn .&\nn .&\nn .&\nn .&\nn .&\nn .&\nn .&\nn .&\nn .&\nn .&\nn 1&\nn .&\nn .&\nn 1&\nn .&\nn 1&\nn .&\nn 1&\nn .&\nn . \\[\sss] 
  15 & 4 &   263424 & .&\nn .&\nn .&\nn .&\nn .&\nn .&\nn .&\nn .&\nn 1&\nn .&\nn 1&\nn .&\nn 1&\nn .&\nn .&\nn .&\nn .&\nn .&\nn 1&\nn .&\nn .&\nn 1&\nn 1&\nn .&\nn .&\nn .&\nn . \\[\sss] 
  16 & 4 &   444528 & .&\nn .&\nn .&\nn .&\nn .&\nn .&\nn .&\nn .&\nn .&\nn .&\nn .&\nn .&\nn .&\nn .&\nn 1&\nn .&\nn 1&\nn .&\nn .&\nn .&\nn 1&\nn .&\nn 1&\nn .&\nn 1&\nn .&\nn . \\[\sss] 
  17 & 4 &   592704 & .&\nn .&\nn .&\nn .&\nn .&\nn .&\nn .&\nn .&\nn .&\nn .&\nn .&\nn .&\nn .&\nn .&\nn 1&\nn .&\nn 1&\nn .&\nn 1&\nn .&\nn .&\nn .&\nn 1&\nn .&\nn .&\nn 1&\nn 1 \\[\sss] 
  18 & 4 &  1185408 & .&\nn .&\nn .&\nn .&\nn .&\nn .&\nn .&\nn .&\nn 1&\nn .&\nn .&\nn .&\nn .&\nn .&\nn 1&\nn .&\nn 1&\nn .&\nn 1&\nn .&\nn .&\nn .&\nn .&\nn .&\nn .&\nn .&\nn . \\[\sss] 
  19 & 4 &  1778112 & .&\nn .&\nn .&\nn .&\nn .&\nn .&\nn .&\nn .&\nn .&\nn .&\nn .&\nn .&\nn .&\nn .&\nn 1&\nn .&\nn 1&\nn .&\nn .&\nn .&\nn 1&\nn 1&\nn .&\nn .&\nn .&\nn .&\nn . \\[\sss] 
  20 & 4 &  1778112 & .&\nn .&\nn .&\nn .&\nn .&\nn .&\nn .&\nn .&\nn 1&\nn .&\nn .&\nn .&\nn .&\nn .&\nn .&\nn .&\nn 1&\nn .&\nn .&\nn .&\nn 1&\nn 1&\nn .&\nn .&\nn .&\nn .&\nn . \\[\sss] 
  21 & 4 &  1778112 & .&\nn .&\nn .&\nn .&\nn .&\nn .&\nn .&\nn .&\nn 1&\nn .&\nn .&\nn .&\nn .&\nn 1&\nn .&\nn .&\nn .&\nn .&\nn .&\nn 1&\nn 1&\nn 1&\nn .&\nn .&\nn 1&\nn .&\nn . \\[\sss] 
  22 & 4 &  2370816 & .&\nn .&\nn .&\nn .&\nn .&\nn .&\nn .&\nn .&\nn 1&\nn .&\nn 1&\nn .&\nn 1&\nn .&\nn .&\nn .&\nn .&\nn .&\nn 1&\nn .&\nn .&\nn 1&\nn 1&\nn .&\nn .&\nn 1&\nn . \\[\sss] 
  23 & 4 &  2370816 & .&\nn .&\nn .&\nn .&\nn .&\nn .&\nn .&\nn .&\nn 1&\nn .&\nn 1&\nn .&\nn 1&\nn .&\nn .&\nn .&\nn .&\nn .&\nn 1&\nn .&\nn .&\nn 1&\nn 1&\nn 1&\nn .&\nn 1&\nn . \\[\sss] 
  24 & 4 &  3556224 & .&\nn .&\nn .&\nn .&\nn .&\nn .&\nn .&\nn .&\nn 1&\nn .&\nn .&\nn .&\nn .&\nn .&\nn 1&\nn .&\nn 1&\nn .&\nn 1&\nn .&\nn .&\nn .&\nn .&\nn .&\nn .&\nn 1&\nn . \\[\sss] 
  25 & 4 &  4741632 & .&\nn .&\nn .&\nn .&\nn .&\nn .&\nn .&\nn .&\nn 1&\nn .&\nn .&\nn .&\nn .&\nn 1&\nn .&\nn .&\nn .&\nn .&\nn 1&\nn .&\nn .&\nn .&\nn .&\nn .&\nn .&\nn 1&\nn . \\[\sss] 
  26 & 4 &  4741632 & .&\nn .&\nn .&\nn .&\nn .&\nn .&\nn .&\nn .&\nn 1&\nn .&\nn .&\nn .&\nn .&\nn 1&\nn .&\nn 1&\nn .&\nn .&\nn 1&\nn .&\nn .&\nn .&\nn .&\nn 1&\nn .&\nn .&\nn . \\[\sss] 
  27 & 4 &  7112448 & .&\nn .&\nn .&\nn .&\nn .&\nn .&\nn .&\nn .&\nn 1&\nn .&\nn .&\nn .&\nn .&\nn .&\nn 1&\nn .&\nn 1&\nn .&\nn 1&\nn .&\nn .&\nn .&\nn 1&\nn .&\nn .&\nn .&\nn . \\[\sss] 
  28 & 4 &  7112448 & .&\nn .&\nn .&\nn .&\nn .&\nn .&\nn .&\nn .&\nn 1&\nn .&\nn .&\nn .&\nn .&\nn 1&\nn .&\nn .&\nn .&\nn .&\nn 1&\nn .&\nn .&\nn .&\nn .&\nn 1&\nn .&\nn 1&\nn . \\[\sss] 
  29 & 4 &  7112448 & .&\nn .&\nn .&\nn .&\nn .&\nn .&\nn .&\nn .&\nn 1&\nn .&\nn .&\nn .&\nn .&\nn 1&\nn .&\nn 1&\nn .&\nn .&\nn .&\nn 1&\nn 1&\nn 1&\nn .&\nn .&\nn .&\nn .&\nn . \\ \midrule 
  30 & 5 &    28224 & .&\nn .&\nn .&\nn .&\nn .&\nn 1&\nn .&\nn 1&\nn .&\nn .&\nn .&\nn 1&\nn .&\nn .&\nn .&\nn 1&\nn .&\nn .&\nn .&\nn 1&\nn .&\nn 1&\nn .&\nn .&\nn .&\nn .&\nn . \\[\sss] 
  31 & 5 &   148176 & .&\nn .&\nn .&\nn .&\nn .&\nn .&\nn .&\nn .&\nn 1&\nn .&\nn .&\nn .&\nn .&\nn .&\nn .&\nn .&\nn 1&\nn .&\nn .&\nn .&\nn 1&\nn .&\nn 1&\nn .&\nn 1&\nn .&\nn . \\[\sss] 
  32 & 5 &   148176 & .&\nn .&\nn .&\nn .&\nn .&\nn .&\nn .&\nn .&\nn 1&\nn .&\nn .&\nn .&\nn .&\nn .&\nn 1&\nn .&\nn 1&\nn .&\nn .&\nn .&\nn 1&\nn .&\nn 1&\nn .&\nn 1&\nn .&\nn . \\[\sss] 
  33 & 5 &   169344 & .&\nn .&\nn .&\nn .&\nn .&\nn .&\nn .&\nn .&\nn .&\nn .&\nn .&\nn 1&\nn .&\nn 1&\nn .&\nn 1&\nn .&\nn .&\nn .&\nn 1&\nn .&\nn 1&\nn .&\nn .&\nn .&\nn 1&\nn 1 \\[\sss] 
  34 & 5 &   592704 & .&\nn .&\nn .&\nn .&\nn .&\nn 1&\nn .&\nn 1&\nn .&\nn .&\nn .&\nn 1&\nn .&\nn .&\nn .&\nn 1&\nn .&\nn .&\nn .&\nn 1&\nn .&\nn 1&\nn 1&\nn .&\nn .&\nn 1&\nn 1 \\[\sss] 
  35 & 5 &  1185408 & .&\nn .&\nn .&\nn .&\nn .&\nn 1&\nn .&\nn 1&\nn .&\nn .&\nn .&\nn 1&\nn .&\nn .&\nn .&\nn 1&\nn .&\nn .&\nn .&\nn 1&\nn .&\nn 1&\nn .&\nn .&\nn .&\nn 1&\nn . \\[\sss] 
  36 & 5 &  1580544 & .&\nn .&\nn .&\nn .&\nn .&\nn 1&\nn .&\nn 1&\nn .&\nn .&\nn .&\nn 1&\nn .&\nn .&\nn .&\nn 1&\nn .&\nn .&\nn 1&\nn .&\nn .&\nn 1&\nn 1&\nn .&\nn .&\nn .&\nn . \\[\sss] 
  37 & 5 &  1580544 & .&\nn .&\nn .&\nn .&\nn .&\nn 1&\nn .&\nn 1&\nn .&\nn .&\nn .&\nn 1&\nn .&\nn .&\nn .&\nn 1&\nn .&\nn .&\nn 1&\nn .&\nn .&\nn 1&\nn 1&\nn .&\nn .&\nn .&\nn 1 \\[\sss] 
  38 & 5 &  1778112 & .&\nn .&\nn .&\nn .&\nn .&\nn .&\nn .&\nn .&\nn 1&\nn .&\nn .&\nn .&\nn .&\nn .&\nn 1&\nn .&\nn 1&\nn .&\nn .&\nn .&\nn 1&\nn 1&\nn .&\nn .&\nn .&\nn .&\nn . \\[\sss] 
  39 & 5 &  1778112 & .&\nn .&\nn .&\nn .&\nn .&\nn .&\nn .&\nn .&\nn 1&\nn .&\nn .&\nn .&\nn .&\nn 1&\nn .&\nn 1&\nn .&\nn .&\nn .&\nn .&\nn 1&\nn 1&\nn .&\nn .&\nn 1&\nn 1&\nn . \\[\sss] 
  40 & 5 &  2370816 & .&\nn .&\nn .&\nn .&\nn .&\nn 1&\nn .&\nn 1&\nn .&\nn .&\nn .&\nn 1&\nn .&\nn .&\nn .&\nn 1&\nn .&\nn .&\nn .&\nn 1&\nn .&\nn 1&\nn 1&\nn .&\nn 1&\nn .&\nn . \\[\sss] 
  41 & 5 &  2370816 & .&\nn .&\nn .&\nn .&\nn .&\nn 1&\nn .&\nn 1&\nn .&\nn .&\nn .&\nn 1&\nn .&\nn 1&\nn .&\nn 1&\nn .&\nn .&\nn 1&\nn .&\nn .&\nn 1&\nn .&\nn .&\nn .&\nn .&\nn 1 \\[\sss] 
  42 & 5 &  2370816 & .&\nn .&\nn .&\nn .&\nn .&\nn 1&\nn .&\nn 1&\nn .&\nn .&\nn .&\nn 1&\nn .&\nn 1&\nn .&\nn 1&\nn .&\nn .&\nn 1&\nn .&\nn .&\nn 1&\nn .&\nn .&\nn 1&\nn .&\nn 1 \\[\sss] 
  43 & 5 &  3556224 & .&\nn .&\nn .&\nn .&\nn .&\nn .&\nn .&\nn .&\nn 1&\nn .&\nn .&\nn .&\nn .&\nn .&\nn 1&\nn .&\nn 1&\nn .&\nn .&\nn 1&\nn .&\nn 1&\nn .&\nn .&\nn .&\nn .&\nn . \\[\sss] 
  44 & 5 &  4741632 & .&\nn .&\nn .&\nn .&\nn .&\nn .&\nn .&\nn .&\nn 1&\nn .&\nn 1&\nn .&\nn 1&\nn .&\nn .&\nn .&\nn .&\nn .&\nn 1&\nn .&\nn .&\nn 1&\nn 1&\nn 1&\nn 1&\nn .&\nn . \\[\sss] 
  45 & 5 &  4741632 & .&\nn .&\nn .&\nn .&\nn .&\nn 1&\nn .&\nn 1&\nn .&\nn .&\nn .&\nn 1&\nn .&\nn 1&\nn .&\nn 1&\nn .&\nn .&\nn .&\nn 1&\nn .&\nn 1&\nn .&\nn .&\nn 1&\nn .&\nn . \\[\sss] 
  46 & 5 &  4741632 & .&\nn .&\nn .&\nn .&\nn .&\nn 1&\nn .&\nn 1&\nn .&\nn .&\nn .&\nn 1&\nn .&\nn 1&\nn .&\nn 1&\nn .&\nn .&\nn 1&\nn .&\nn .&\nn .&\nn .&\nn .&\nn .&\nn 1&\nn . \\[\sss] 
  47 & 5 &  4741632 & .&\nn .&\nn .&\nn .&\nn .&\nn 1&\nn .&\nn 1&\nn .&\nn .&\nn .&\nn 1&\nn .&\nn 1&\nn .&\nn 1&\nn .&\nn .&\nn 1&\nn .&\nn .&\nn .&\nn .&\nn .&\nn .&\nn 1&\nn 1 \\[\sss] 
  48 & 5 &  4741632 & .&\nn .&\nn .&\nn .&\nn .&\nn 1&\nn .&\nn 1&\nn .&\nn .&\nn .&\nn 1&\nn .&\nn 1&\nn .&\nn 1&\nn .&\nn .&\nn 1&\nn .&\nn .&\nn 1&\nn .&\nn .&\nn .&\nn 1&\nn . \\[\sss] 
  49 & 5 &  4741632 & .&\nn .&\nn .&\nn .&\nn .&\nn 1&\nn .&\nn 1&\nn .&\nn .&\nn .&\nn 1&\nn .&\nn 1&\nn .&\nn 1&\nn .&\nn .&\nn 1&\nn .&\nn .&\nn 1&\nn .&\nn .&\nn 1&\nn 1&\nn . \\[\sss] 
  50 & 5 &  7112448 & .&\nn .&\nn .&\nn .&\nn .&\nn .&\nn .&\nn .&\nn 1&\nn .&\nn .&\nn 1&\nn .&\nn 1&\nn .&\nn 1&\nn .&\nn .&\nn 1&\nn .&\nn .&\nn .&\nn .&\nn 1&\nn .&\nn 1&\nn . \\[\sss] 
  51 & 5 & 14224896 & .&\nn .&\nn .&\nn .&\nn .&\nn .&\nn .&\nn .&\nn 1&\nn .&\nn .&\nn .&\nn .&\nn 1&\nn .&\nn 1&\nn .&\nn .&\nn 1&\nn .&\nn .&\nn .&\nn .&\nn 1&\nn .&\nn 1&\nn . \\[\sss] 
  52 & 5 & 14224896 & .&\nn .&\nn .&\nn .&\nn .&\nn .&\nn .&\nn .&\nn 1&\nn .&\nn .&\nn 1&\nn .&\nn 1&\nn .&\nn 1&\nn .&\nn .&\nn .&\nn 1&\nn .&\nn 1&\nn .&\nn .&\nn .&\nn .&\nn . \\ \midrule 
  53 & 6 &    32256 & .&\nn .&\nn 1&\nn .&\nn 1&\nn .&\nn 1&\nn .&\nn .&\nn .&\nn 1&\nn .&\nn 1&\nn .&\nn .&\nn .&\nn 1&\nn 1&\nn 1&\nn .&\nn .&\nn .&\nn 1&\nn 1&\nn 1&\nn 1&\nn . \\[\sss] 
  54 & 6 &   197568 & .&\nn .&\nn .&\nn .&\nn .&\nn 1&\nn .&\nn 1&\nn .&\nn .&\nn .&\nn 1&\nn .&\nn .&\nn .&\nn 1&\nn .&\nn .&\nn .&\nn 1&\nn .&\nn 1&\nn .&\nn .&\nn .&\nn .&\nn 1 \\[\sss] 
  55 & 6 &   395136 & .&\nn .&\nn .&\nn .&\nn .&\nn 1&\nn .&\nn 1&\nn .&\nn .&\nn .&\nn 1&\nn .&\nn 1&\nn .&\nn 1&\nn .&\nn .&\nn .&\nn 1&\nn .&\nn 1&\nn .&\nn .&\nn .&\nn 1&\nn 1 \\ \midrule
\end{tabular}
\caption{Large orbits of $3 \times 3 \times 3$ tensors}
\label{table333}
\end{table}  
   
%%%%%%%%%%%%%%%%%%%%%%%%%%%%%%%%%%%%%%%%%%%%%%%%%%%%%%%%%%%%%%%%%%%%%%%%%%%%%%

\section{Formats $p \times 2 \times 2 \times 2$ ($p = 2, 3$)}

For $p = 2$ the classification appears in \cite{BS}; there are $2^{16} = 65536$ tensors.  
The small group order is $6^4 = 1296$, and the large group order is $24 \cdot 6^4 = 31104$
(we include all permutations of the four equal directions).

For $p = 3$, the large group acting on the tensors is
  \[
  G = \big( 
  GL(3,\mathbb{F}_2) \times 
  GL(2,\mathbb{F}_2) \times 
  GL(2,\mathbb{F}_2) \times 
  GL(2,\mathbb{F}_2) 
  \big) 
  \rtimes S_3,
  \]
where $S_3$ permutes the last three directions.
We follow the same approach as the previous section,
but now there are $2^{24} = 16777216$ tensors.
The small group order is $168 \cdot 6^3 = 36288$, 
and the large group order is $6 \cdot 168 \cdot 6^3 = 217728$.
There are 696 (nonzero) small orbits and 212 (nonzero) large orbits, with maximum rank 6:
\[
\begin{array}{lrrrrrrr}
\text{rank} 
& 0 & 1 & 2 & 3 & 4 & 5 & 6 
\\
\text{$\#$ small} 
& 1 & 1 & 11 & 49 & 234 & 367 & 34
\\
\text{$\#$ large} 
& 1 & 1 & 5 & 21 & 72 & 100 & 13
\\
\text{$\#$ tensors} 
& 1 & 189 & 13608 & 434028 & 5143446 & 10460016 & 725928
\\
\text{percent}
& 0.0000
& 0.0011
& 0.0811
& 2.5870
& 30.6573
& 62.3466
& 4.3269
\end{array}
\]
For the large orbit sizes and canonical forms, see Tables \ref{table3222part1} and \ref{table3222part2}.
This computation took less than 40 minutes.

\renewcommand{\nn}{\!\!\!\!}

\renewcommand{\sss}{-2pt}

\begin{table} \tiny
\begin{tabular}{rrrcccccccccccccccccccccccc}
$\#$ & rank & size & \multicolumn{24}{c}{canonical form} \\ \midrule
    1 & 1 &    189 & .&\nn .&\nn .&\nn .&\nn .&\nn .&\nn .&\nn .&\nn .&\nn .&\nn .&\nn .&\nn .&\nn .&\nn .&\nn .&\nn .&\nn .&\nn .&\nn .&\nn .&\nn .&\nn .&\nn 1 \\ \midrule 
    2 & 2 &    378 & .&\nn .&\nn .&\nn .&\nn .&\nn .&\nn .&\nn .&\nn .&\nn .&\nn .&\nn .&\nn .&\nn .&\nn .&\nn .&\nn .&\nn .&\nn .&\nn .&\nn .&\nn 1&\nn 1&\nn . \\[\sss] 
    3 & 2 &    756 & .&\nn .&\nn .&\nn .&\nn .&\nn .&\nn .&\nn .&\nn .&\nn .&\nn .&\nn .&\nn .&\nn .&\nn .&\nn .&\nn .&\nn .&\nn .&\nn 1&\nn 1&\nn .&\nn .&\nn . \\[\sss] 
    4 & 2 &   1134 & .&\nn .&\nn .&\nn .&\nn .&\nn .&\nn .&\nn .&\nn .&\nn .&\nn .&\nn .&\nn .&\nn .&\nn .&\nn 1&\nn .&\nn .&\nn .&\nn .&\nn .&\nn .&\nn 1&\nn . \\[\sss] 
    5 & 2 &   4536 & .&\nn .&\nn .&\nn .&\nn .&\nn .&\nn .&\nn .&\nn .&\nn .&\nn .&\nn .&\nn .&\nn .&\nn .&\nn 1&\nn 1&\nn .&\nn .&\nn .&\nn .&\nn .&\nn .&\nn . \\[\sss] 
    6 & 2 &   6804 & .&\nn .&\nn .&\nn .&\nn .&\nn .&\nn .&\nn .&\nn .&\nn .&\nn .&\nn .&\nn .&\nn .&\nn .&\nn 1&\nn .&\nn .&\nn .&\nn .&\nn 1&\nn .&\nn .&\nn . \\ \midrule 
    7 & 3 &     84 & .&\nn .&\nn .&\nn .&\nn .&\nn .&\nn .&\nn .&\nn .&\nn .&\nn .&\nn .&\nn .&\nn .&\nn .&\nn .&\nn .&\nn 1&\nn 1&\nn .&\nn 1&\nn .&\nn 1&\nn 1 \\[\sss] 
    8 & 3 &    378 & .&\nn .&\nn .&\nn .&\nn .&\nn .&\nn .&\nn .&\nn .&\nn .&\nn .&\nn .&\nn .&\nn .&\nn .&\nn .&\nn .&\nn .&\nn .&\nn 1&\nn .&\nn 1&\nn 1&\nn . \\[\sss] 
    9 & 3 &    756 & .&\nn .&\nn .&\nn .&\nn .&\nn .&\nn .&\nn .&\nn .&\nn .&\nn .&\nn .&\nn .&\nn 1&\nn 1&\nn .&\nn .&\nn .&\nn .&\nn .&\nn 1&\nn .&\nn 1&\nn 1 \\[\sss] 
   10 & 3 &   1512 & .&\nn .&\nn .&\nn .&\nn .&\nn .&\nn .&\nn .&\nn .&\nn .&\nn .&\nn 1&\nn 1&\nn .&\nn .&\nn .&\nn 1&\nn 1&\nn 1&\nn .&\nn 1&\nn 1&\nn 1&\nn 1 \\[\sss] 
   11 & 3 &   3402 & .&\nn .&\nn .&\nn .&\nn .&\nn .&\nn .&\nn .&\nn .&\nn .&\nn .&\nn .&\nn .&\nn .&\nn .&\nn 1&\nn .&\nn .&\nn .&\nn .&\nn .&\nn 1&\nn 1&\nn . \\[\sss] 
   12 & 3 &   4536 & .&\nn .&\nn .&\nn .&\nn .&\nn .&\nn .&\nn .&\nn .&\nn .&\nn .&\nn .&\nn .&\nn .&\nn .&\nn 1&\nn .&\nn 1&\nn 1&\nn .&\nn 1&\nn .&\nn 1&\nn . \\[\sss] 
   13 & 3 &   4536 & .&\nn .&\nn .&\nn .&\nn .&\nn .&\nn .&\nn .&\nn .&\nn .&\nn .&\nn .&\nn .&\nn 1&\nn 1&\nn .&\nn .&\nn .&\nn .&\nn 1&\nn .&\nn .&\nn 1&\nn . \\[\sss] 
   14 & 3 &   6048 & .&\nn .&\nn .&\nn .&\nn .&\nn .&\nn .&\nn 1&\nn .&\nn 1&\nn 1&\nn .&\nn 1&\nn .&\nn 1&\nn .&\nn 1&\nn .&\nn 1&\nn .&\nn 1&\nn .&\nn 1&\nn . \\[\sss] 
   15 & 3 &   9072 & .&\nn .&\nn .&\nn .&\nn .&\nn .&\nn .&\nn 1&\nn .&\nn .&\nn .&\nn .&\nn .&\nn 1&\nn 1&\nn .&\nn .&\nn .&\nn .&\nn .&\nn 1&\nn .&\nn .&\nn . \\[\sss] 
   16 & 3 &  13608 & .&\nn .&\nn .&\nn .&\nn .&\nn .&\nn .&\nn .&\nn .&\nn .&\nn .&\nn .&\nn .&\nn .&\nn .&\nn 1&\nn .&\nn .&\nn .&\nn 1&\nn 1&\nn .&\nn .&\nn . \\[\sss] 
   17 & 3 &  13608 & .&\nn .&\nn .&\nn .&\nn .&\nn .&\nn .&\nn .&\nn .&\nn .&\nn .&\nn .&\nn .&\nn .&\nn .&\nn 1&\nn .&\nn .&\nn 1&\nn .&\nn 1&\nn .&\nn .&\nn . \\[\sss] 
   18 & 3 &  13608 & .&\nn .&\nn .&\nn .&\nn .&\nn .&\nn .&\nn .&\nn .&\nn .&\nn .&\nn .&\nn .&\nn .&\nn .&\nn 1&\nn .&\nn .&\nn 1&\nn .&\nn 1&\nn 1&\nn .&\nn . \\[\sss] 
   19 & 3 &  13608 & .&\nn .&\nn .&\nn .&\nn .&\nn .&\nn .&\nn .&\nn .&\nn .&\nn .&\nn .&\nn .&\nn 1&\nn 1&\nn .&\nn .&\nn .&\nn .&\nn 1&\nn 1&\nn .&\nn 1&\nn . \\[\sss] 
   20 & 3 &  13608 & .&\nn .&\nn .&\nn .&\nn .&\nn .&\nn .&\nn 1&\nn .&\nn .&\nn .&\nn .&\nn .&\nn .&\nn 1&\nn .&\nn .&\nn .&\nn .&\nn .&\nn .&\nn 1&\nn .&\nn . \\[\sss] 
   21 & 3 &  27216 & .&\nn .&\nn .&\nn .&\nn .&\nn .&\nn .&\nn .&\nn .&\nn .&\nn .&\nn .&\nn .&\nn .&\nn .&\nn 1&\nn 1&\nn .&\nn .&\nn .&\nn .&\nn .&\nn 1&\nn . \\[\sss] 
   22 & 3 &  27216 & .&\nn .&\nn .&\nn .&\nn .&\nn .&\nn .&\nn .&\nn .&\nn .&\nn .&\nn .&\nn .&\nn .&\nn .&\nn 1&\nn 1&\nn .&\nn .&\nn .&\nn .&\nn 1&\nn 1&\nn . \\[\sss] 
   23 & 3 &  27216 & .&\nn .&\nn .&\nn .&\nn .&\nn .&\nn .&\nn .&\nn .&\nn .&\nn .&\nn .&\nn .&\nn 1&\nn 1&\nn .&\nn .&\nn .&\nn 1&\nn .&\nn .&\nn .&\nn .&\nn 1 \\[\sss] 
   24 & 3 &  36288 & .&\nn .&\nn .&\nn .&\nn .&\nn .&\nn .&\nn 1&\nn .&\nn .&\nn .&\nn .&\nn 1&\nn .&\nn .&\nn .&\nn .&\nn .&\nn 1&\nn .&\nn .&\nn .&\nn .&\nn . \\[\sss] 
   25 & 3 &  54432 & .&\nn .&\nn .&\nn .&\nn .&\nn .&\nn .&\nn 1&\nn .&\nn .&\nn .&\nn .&\nn .&\nn .&\nn 1&\nn .&\nn .&\nn 1&\nn .&\nn .&\nn .&\nn .&\nn .&\nn . \\[\sss] 
   26 & 3 &  54432 & .&\nn .&\nn .&\nn .&\nn .&\nn .&\nn .&\nn 1&\nn .&\nn .&\nn .&\nn .&\nn 1&\nn .&\nn .&\nn .&\nn 1&\nn 1&\nn 1&\nn 1&\nn .&\nn .&\nn .&\nn . \\[\sss] 
   27 & 3 & 108864 & .&\nn .&\nn .&\nn .&\nn .&\nn .&\nn .&\nn 1&\nn .&\nn .&\nn .&\nn .&\nn 1&\nn .&\nn .&\nn .&\nn .&\nn .&\nn 1&\nn 1&\nn .&\nn .&\nn .&\nn . \\ \midrule 
   28 & 4 &    756 & .&\nn .&\nn .&\nn .&\nn .&\nn .&\nn .&\nn .&\nn .&\nn .&\nn .&\nn .&\nn .&\nn 1&\nn 1&\nn .&\nn .&\nn 1&\nn 1&\nn .&\nn .&\nn .&\nn .&\nn . \\[\sss] 
   29 & 4 &   1134 & .&\nn .&\nn .&\nn .&\nn .&\nn .&\nn .&\nn .&\nn .&\nn .&\nn .&\nn .&\nn .&\nn .&\nn .&\nn 1&\nn .&\nn .&\nn .&\nn 1&\nn .&\nn 1&\nn 1&\nn . \\[\sss] 
   30 & 4 &   4536 & .&\nn .&\nn .&\nn .&\nn .&\nn .&\nn .&\nn .&\nn .&\nn .&\nn .&\nn .&\nn .&\nn .&\nn .&\nn 1&\nn .&\nn 1&\nn 1&\nn .&\nn 1&\nn .&\nn .&\nn . \\[\sss] 
   31 & 4 &   4536 & .&\nn .&\nn .&\nn .&\nn .&\nn .&\nn .&\nn .&\nn .&\nn .&\nn .&\nn .&\nn .&\nn 1&\nn 1&\nn .&\nn 1&\nn .&\nn 1&\nn 1&\nn .&\nn .&\nn .&\nn . \\[\sss] 
   32 & 4 &   4536 & .&\nn .&\nn .&\nn .&\nn .&\nn .&\nn .&\nn .&\nn .&\nn .&\nn .&\nn 1&\nn 1&\nn .&\nn .&\nn .&\nn .&\nn .&\nn 1&\nn .&\nn 1&\nn 1&\nn .&\nn . \\[\sss] 
   33 & 4 &   4536 & .&\nn .&\nn .&\nn .&\nn .&\nn .&\nn .&\nn 1&\nn .&\nn .&\nn .&\nn .&\nn .&\nn 1&\nn 1&\nn .&\nn .&\nn .&\nn .&\nn 1&\nn .&\nn .&\nn 1&\nn . \\[\sss] 
   34 & 4 &   6804 & .&\nn .&\nn .&\nn .&\nn .&\nn .&\nn .&\nn .&\nn .&\nn .&\nn .&\nn .&\nn .&\nn 1&\nn 1&\nn .&\nn .&\nn 1&\nn 1&\nn .&\nn .&\nn .&\nn .&\nn 1 \\[\sss] 
   35 & 4 &   6804 & .&\nn .&\nn .&\nn .&\nn .&\nn .&\nn .&\nn .&\nn .&\nn .&\nn .&\nn .&\nn .&\nn 1&\nn 1&\nn .&\nn .&\nn 1&\nn 1&\nn .&\nn .&\nn .&\nn 1&\nn . \\[\sss] 
   36 & 4 &   6804 & .&\nn .&\nn .&\nn .&\nn .&\nn .&\nn .&\nn .&\nn .&\nn .&\nn .&\nn .&\nn .&\nn 1&\nn 1&\nn .&\nn .&\nn 1&\nn 1&\nn 1&\nn .&\nn .&\nn .&\nn . \\[\sss] 
   37 & 4 &   9072 & .&\nn .&\nn .&\nn .&\nn .&\nn .&\nn .&\nn .&\nn .&\nn .&\nn .&\nn .&\nn .&\nn .&\nn .&\nn 1&\nn 1&\nn .&\nn .&\nn 1&\nn .&\nn 1&\nn 1&\nn . \\[\sss] 
   38 & 4 &   9072 & .&\nn .&\nn .&\nn .&\nn .&\nn 1&\nn 1&\nn .&\nn .&\nn .&\nn .&\nn 1&\nn .&\nn .&\nn 1&\nn .&\nn 1&\nn .&\nn .&\nn .&\nn .&\nn .&\nn 1&\nn . \\[\sss] 
   39 & 4 &  13608 & .&\nn .&\nn .&\nn .&\nn .&\nn .&\nn .&\nn .&\nn .&\nn .&\nn .&\nn .&\nn .&\nn 1&\nn 1&\nn .&\nn .&\nn .&\nn .&\nn 1&\nn 1&\nn .&\nn .&\nn . \\[\sss] 
   40 & 4 &  13608 & .&\nn .&\nn .&\nn .&\nn .&\nn .&\nn .&\nn .&\nn .&\nn .&\nn .&\nn .&\nn .&\nn 1&\nn 1&\nn .&\nn .&\nn 1&\nn 1&\nn 1&\nn .&\nn .&\nn 1&\nn . \\[\sss] 
   41 & 4 &  13608 & .&\nn .&\nn .&\nn .&\nn .&\nn .&\nn .&\nn .&\nn .&\nn .&\nn .&\nn .&\nn .&\nn 1&\nn 1&\nn .&\nn .&\nn 1&\nn 1&\nn 1&\nn 1&\nn .&\nn .&\nn . \\[\sss] 
   42 & 4 &  13608 & .&\nn .&\nn .&\nn .&\nn .&\nn .&\nn .&\nn .&\nn .&\nn .&\nn .&\nn .&\nn .&\nn 1&\nn 1&\nn .&\nn .&\nn 1&\nn 1&\nn 1&\nn 1&\nn .&\nn 1&\nn . \\[\sss] 
   43 & 4 &  13608 & .&\nn .&\nn .&\nn .&\nn .&\nn .&\nn .&\nn .&\nn .&\nn .&\nn .&\nn 1&\nn .&\nn 1&\nn 1&\nn .&\nn 1&\nn .&\nn .&\nn .&\nn 1&\nn .&\nn 1&\nn 1 \\[\sss] 
   44 & 4 &  13608 & .&\nn .&\nn .&\nn .&\nn .&\nn .&\nn .&\nn .&\nn .&\nn .&\nn .&\nn 1&\nn 1&\nn .&\nn .&\nn .&\nn .&\nn 1&\nn 1&\nn .&\nn 1&\nn .&\nn 1&\nn 1 \\[\sss] 
   45 & 4 &  13608 & .&\nn .&\nn .&\nn .&\nn .&\nn .&\nn .&\nn 1&\nn .&\nn .&\nn .&\nn .&\nn .&\nn .&\nn 1&\nn .&\nn .&\nn .&\nn .&\nn 1&\nn .&\nn 1&\nn .&\nn . \\[\sss] 
   46 & 4 &  27216 & .&\nn .&\nn .&\nn .&\nn .&\nn .&\nn .&\nn .&\nn .&\nn .&\nn .&\nn .&\nn .&\nn 1&\nn 1&\nn .&\nn .&\nn .&\nn 1&\nn .&\nn 1&\nn .&\nn .&\nn 1 \\[\sss] 
   47 & 4 &  27216 & .&\nn .&\nn .&\nn .&\nn .&\nn .&\nn .&\nn .&\nn .&\nn .&\nn .&\nn .&\nn .&\nn 1&\nn 1&\nn .&\nn 1&\nn .&\nn 1&\nn 1&\nn .&\nn .&\nn .&\nn 1 \\[\sss] 
   48 & 4 &  27216 & .&\nn .&\nn .&\nn .&\nn .&\nn .&\nn .&\nn .&\nn .&\nn .&\nn .&\nn 1&\nn .&\nn 1&\nn 1&\nn .&\nn .&\nn .&\nn 1&\nn .&\nn 1&\nn 1&\nn .&\nn . \\[\sss] 
   49 & 4 &  27216 & .&\nn .&\nn .&\nn .&\nn .&\nn .&\nn .&\nn .&\nn .&\nn .&\nn .&\nn 1&\nn .&\nn 1&\nn 1&\nn .&\nn 1&\nn .&\nn .&\nn .&\nn .&\nn .&\nn 1&\nn 1 \\[\sss] 
   50 & 4 &  27216 & .&\nn .&\nn .&\nn .&\nn .&\nn .&\nn .&\nn .&\nn .&\nn .&\nn .&\nn 1&\nn 1&\nn .&\nn .&\nn .&\nn .&\nn .&\nn 1&\nn .&\nn 1&\nn 1&\nn .&\nn 1 \\[\sss] 
   51 & 4 &  27216 & .&\nn .&\nn .&\nn .&\nn .&\nn .&\nn .&\nn 1&\nn .&\nn .&\nn .&\nn .&\nn .&\nn .&\nn 1&\nn .&\nn .&\nn .&\nn .&\nn 1&\nn 1&\nn .&\nn .&\nn . \\[\sss] 
   52 & 4 &  27216 & .&\nn .&\nn .&\nn .&\nn .&\nn .&\nn .&\nn 1&\nn .&\nn .&\nn .&\nn .&\nn .&\nn 1&\nn 1&\nn .&\nn .&\nn .&\nn .&\nn 1&\nn 1&\nn .&\nn .&\nn . \\[\sss] 
   53 & 4 &  27216 & .&\nn .&\nn .&\nn .&\nn .&\nn .&\nn .&\nn 1&\nn .&\nn .&\nn .&\nn .&\nn 1&\nn .&\nn .&\nn .&\nn 1&\nn .&\nn .&\nn 1&\nn .&\nn .&\nn .&\nn . \\[\sss] 
   54 & 4 &  36288 & .&\nn .&\nn .&\nn .&\nn .&\nn .&\nn .&\nn 1&\nn .&\nn .&\nn 1&\nn .&\nn 1&\nn .&\nn .&\nn .&\nn .&\nn 1&\nn .&\nn .&\nn 1&\nn .&\nn .&\nn . \\[\sss] 
   55 & 4 &  36288 & .&\nn .&\nn .&\nn .&\nn .&\nn .&\nn .&\nn 1&\nn .&\nn .&\nn 1&\nn .&\nn 1&\nn 1&\nn .&\nn .&\nn .&\nn 1&\nn .&\nn .&\nn 1&\nn 1&\nn 1&\nn . \\[\sss] 
   56 & 4 &  36288 & .&\nn .&\nn .&\nn .&\nn .&\nn 1&\nn 1&\nn .&\nn .&\nn .&\nn .&\nn 1&\nn .&\nn .&\nn 1&\nn .&\nn 1&\nn .&\nn .&\nn .&\nn .&\nn .&\nn 1&\nn 1 \\[\sss] 
   57 & 4 &  36288 & .&\nn .&\nn .&\nn .&\nn .&\nn 1&\nn 1&\nn .&\nn .&\nn .&\nn 1&\nn .&\nn .&\nn .&\nn .&\nn 1&\nn 1&\nn 1&\nn .&\nn .&\nn .&\nn .&\nn 1&\nn 1 \\[\sss] 
   58 & 4 &  54432 & .&\nn .&\nn .&\nn .&\nn .&\nn .&\nn .&\nn 1&\nn .&\nn .&\nn .&\nn .&\nn .&\nn .&\nn 1&\nn .&\nn .&\nn 1&\nn 1&\nn .&\nn 1&\nn .&\nn .&\nn . \\[\sss] 
   59 & 4 &  54432 & .&\nn .&\nn .&\nn .&\nn .&\nn .&\nn .&\nn 1&\nn .&\nn .&\nn .&\nn .&\nn 1&\nn .&\nn .&\nn .&\nn .&\nn .&\nn .&\nn 1&\nn .&\nn 1&\nn 1&\nn . \\[\sss] 
   60 & 4 &  54432 & .&\nn .&\nn .&\nn .&\nn .&\nn .&\nn .&\nn 1&\nn .&\nn .&\nn .&\nn .&\nn 1&\nn .&\nn .&\nn .&\nn .&\nn 1&\nn 1&\nn .&\nn .&\nn .&\nn .&\nn . \\[\sss] 
   61 & 4 &  54432 & .&\nn .&\nn .&\nn .&\nn .&\nn .&\nn .&\nn 1&\nn .&\nn .&\nn .&\nn .&\nn 1&\nn .&\nn .&\nn .&\nn .&\nn 1&\nn 1&\nn .&\nn .&\nn .&\nn 1&\nn . \\[\sss] 
   62 & 4 &  54432 & .&\nn .&\nn .&\nn .&\nn .&\nn .&\nn .&\nn 1&\nn .&\nn .&\nn .&\nn .&\nn 1&\nn .&\nn .&\nn .&\nn 1&\nn .&\nn .&\nn 1&\nn .&\nn .&\nn 1&\nn . \\[\sss] 
   63 & 4 &  54432 & .&\nn .&\nn .&\nn .&\nn .&\nn .&\nn .&\nn 1&\nn .&\nn .&\nn .&\nn .&\nn 1&\nn .&\nn .&\nn .&\nn 1&\nn 1&\nn 1&\nn 1&\nn .&\nn .&\nn 1&\nn . \\[\sss] 
   64 & 4 &  54432 & .&\nn .&\nn .&\nn .&\nn .&\nn .&\nn .&\nn 1&\nn .&\nn .&\nn .&\nn 1&\nn 1&\nn .&\nn .&\nn .&\nn 1&\nn .&\nn .&\nn .&\nn .&\nn .&\nn .&\nn . \\[\sss] 
   65 & 4 &  54432 & .&\nn .&\nn .&\nn .&\nn .&\nn .&\nn .&\nn 1&\nn .&\nn .&\nn 1&\nn .&\nn 1&\nn .&\nn .&\nn .&\nn 1&\nn .&\nn .&\nn .&\nn 1&\nn .&\nn 1&\nn . \\[\sss] 
   66 & 4 &  54432 & .&\nn .&\nn .&\nn .&\nn .&\nn .&\nn .&\nn 1&\nn .&\nn .&\nn 1&\nn .&\nn 1&\nn .&\nn .&\nn .&\nn 1&\nn 1&\nn .&\nn .&\nn .&\nn .&\nn .&\nn . \\[\sss] 
   67 & 4 &  54432 & .&\nn .&\nn .&\nn .&\nn .&\nn .&\nn .&\nn 1&\nn .&\nn 1&\nn 1&\nn .&\nn 1&\nn .&\nn 1&\nn .&\nn 1&\nn .&\nn 1&\nn .&\nn 1&\nn .&\nn .&\nn . \\[\sss] 
   68 & 4 &  54432 & .&\nn .&\nn .&\nn .&\nn .&\nn 1&\nn 1&\nn .&\nn .&\nn .&\nn .&\nn .&\nn 1&\nn .&\nn 1&\nn 1&\nn .&\nn .&\nn .&\nn 1&\nn .&\nn .&\nn 1&\nn . \\[\sss] 
   69 & 4 &  54432 & .&\nn .&\nn .&\nn .&\nn .&\nn 1&\nn 1&\nn .&\nn .&\nn .&\nn .&\nn 1&\nn .&\nn .&\nn 1&\nn .&\nn .&\nn .&\nn 1&\nn .&\nn 1&\nn .&\nn .&\nn . \\[\sss] 
   70 & 4 &  54432 & .&\nn .&\nn .&\nn .&\nn .&\nn 1&\nn 1&\nn .&\nn .&\nn .&\nn .&\nn 1&\nn .&\nn .&\nn 1&\nn .&\nn 1&\nn .&\nn .&\nn .&\nn 1&\nn .&\nn 1&\nn . \\[\sss] 
   71 & 4 & 108864 & .&\nn .&\nn .&\nn .&\nn .&\nn .&\nn .&\nn 1&\nn .&\nn .&\nn .&\nn .&\nn .&\nn .&\nn 1&\nn .&\nn .&\nn 1&\nn .&\nn .&\nn 1&\nn .&\nn .&\nn . \\[\sss] 
   72 & 4 & 108864 & .&\nn .&\nn .&\nn .&\nn .&\nn .&\nn .&\nn 1&\nn .&\nn .&\nn .&\nn .&\nn .&\nn 1&\nn 1&\nn .&\nn .&\nn .&\nn 1&\nn .&\nn .&\nn .&\nn .&\nn . \\[\sss] 
   73 & 4 & 108864 & .&\nn .&\nn .&\nn .&\nn .&\nn .&\nn .&\nn 1&\nn .&\nn .&\nn .&\nn .&\nn .&\nn 1&\nn 1&\nn .&\nn .&\nn .&\nn 1&\nn .&\nn 1&\nn .&\nn .&\nn . \\[\sss] 
   74 & 4 & 108864 & .&\nn .&\nn .&\nn .&\nn .&\nn .&\nn .&\nn 1&\nn .&\nn .&\nn .&\nn .&\nn .&\nn 1&\nn 1&\nn .&\nn 1&\nn .&\nn .&\nn .&\nn .&\nn .&\nn .&\nn . \\[\sss] 
   75 & 4 & 108864 & .&\nn .&\nn .&\nn .&\nn .&\nn .&\nn .&\nn 1&\nn .&\nn .&\nn .&\nn .&\nn .&\nn 1&\nn 1&\nn .&\nn 1&\nn .&\nn .&\nn .&\nn .&\nn .&\nn 1&\nn . \\[\sss] 
   76 & 4 & 108864 & .&\nn .&\nn .&\nn .&\nn .&\nn .&\nn .&\nn 1&\nn .&\nn .&\nn .&\nn .&\nn 1&\nn .&\nn .&\nn .&\nn .&\nn .&\nn 1&\nn .&\nn .&\nn 1&\nn .&\nn . \\[\sss] 
   77 & 4 & 108864 & .&\nn .&\nn .&\nn .&\nn .&\nn .&\nn .&\nn 1&\nn .&\nn .&\nn .&\nn .&\nn 1&\nn .&\nn .&\nn .&\nn .&\nn 1&\nn 1&\nn 1&\nn .&\nn .&\nn .&\nn . \\[\sss] 
   78 & 4 & 108864 & .&\nn .&\nn .&\nn .&\nn .&\nn .&\nn .&\nn 1&\nn .&\nn .&\nn .&\nn .&\nn 1&\nn .&\nn .&\nn .&\nn .&\nn 1&\nn 1&\nn 1&\nn .&\nn .&\nn 1&\nn . \\[\sss] 
   79 & 4 & 108864 & .&\nn .&\nn .&\nn .&\nn .&\nn .&\nn .&\nn 1&\nn .&\nn .&\nn .&\nn .&\nn 1&\nn .&\nn .&\nn .&\nn 1&\nn .&\nn 1&\nn 1&\nn .&\nn .&\nn .&\nn . \\[\sss] 
   80 & 4 & 108864 & .&\nn .&\nn .&\nn .&\nn .&\nn .&\nn .&\nn 1&\nn .&\nn .&\nn .&\nn .&\nn 1&\nn .&\nn .&\nn .&\nn 1&\nn .&\nn 1&\nn 1&\nn .&\nn 1&\nn .&\nn . \\[\sss] 
   81 & 4 & 108864 & .&\nn .&\nn .&\nn .&\nn .&\nn .&\nn .&\nn 1&\nn .&\nn .&\nn .&\nn 1&\nn 1&\nn .&\nn .&\nn .&\nn 1&\nn 1&\nn 1&\nn .&\nn .&\nn .&\nn .&\nn . \\[\sss] 
   82 & 4 & 108864 & .&\nn .&\nn .&\nn .&\nn .&\nn .&\nn .&\nn 1&\nn .&\nn .&\nn .&\nn 1&\nn 1&\nn .&\nn .&\nn .&\nn 1&\nn 1&\nn 1&\nn .&\nn .&\nn 1&\nn 1&\nn . \\[\sss] 
   83 & 4 & 108864 & .&\nn .&\nn .&\nn .&\nn .&\nn .&\nn .&\nn 1&\nn .&\nn .&\nn 1&\nn .&\nn 1&\nn .&\nn .&\nn .&\nn .&\nn 1&\nn .&\nn .&\nn 1&\nn .&\nn 1&\nn . \\[\sss] 
   84 & 4 & 108864 & .&\nn .&\nn .&\nn .&\nn .&\nn .&\nn .&\nn 1&\nn .&\nn .&\nn 1&\nn .&\nn 1&\nn .&\nn .&\nn .&\nn .&\nn 1&\nn .&\nn .&\nn 1&\nn 1&\nn .&\nn . \\[\sss] 
   85 & 4 & 108864 & .&\nn .&\nn .&\nn .&\nn .&\nn .&\nn .&\nn 1&\nn .&\nn .&\nn 1&\nn .&\nn 1&\nn .&\nn .&\nn .&\nn 1&\nn .&\nn .&\nn .&\nn 1&\nn 1&\nn .&\nn . \\[\sss] 
   86 & 4 & 108864 & .&\nn .&\nn .&\nn .&\nn .&\nn .&\nn .&\nn 1&\nn .&\nn .&\nn 1&\nn .&\nn 1&\nn .&\nn .&\nn .&\nn 1&\nn 1&\nn .&\nn .&\nn 1&\nn .&\nn .&\nn . \\[\sss] 
   87 & 4 & 108864 & .&\nn .&\nn .&\nn .&\nn .&\nn .&\nn .&\nn 1&\nn .&\nn .&\nn 1&\nn .&\nn 1&\nn .&\nn .&\nn .&\nn 1&\nn 1&\nn .&\nn .&\nn 1&\nn .&\nn 1&\nn . \\[\sss] 
   88 & 4 & 108864 & .&\nn .&\nn .&\nn .&\nn .&\nn .&\nn .&\nn 1&\nn .&\nn .&\nn 1&\nn .&\nn 1&\nn .&\nn .&\nn .&\nn 1&\nn 1&\nn .&\nn .&\nn 1&\nn 1&\nn .&\nn . \\[\sss] 
   89 & 4 & 108864 & .&\nn .&\nn .&\nn .&\nn .&\nn .&\nn .&\nn 1&\nn .&\nn .&\nn 1&\nn .&\nn 1&\nn 1&\nn .&\nn .&\nn 1&\nn .&\nn .&\nn .&\nn 1&\nn .&\nn .&\nn . \\[\sss] 
   90 & 4 & 108864 & .&\nn .&\nn .&\nn .&\nn .&\nn .&\nn .&\nn 1&\nn .&\nn 1&\nn 1&\nn .&\nn 1&\nn .&\nn 1&\nn .&\nn 1&\nn .&\nn .&\nn .&\nn .&\nn 1&\nn .&\nn . \\[\sss] 
   91 & 4 & 108864 & .&\nn .&\nn .&\nn .&\nn .&\nn 1&\nn 1&\nn .&\nn .&\nn .&\nn .&\nn 1&\nn .&\nn .&\nn 1&\nn .&\nn 1&\nn .&\nn .&\nn .&\nn 1&\nn .&\nn 1&\nn 1 \\[\sss] 
   92 & 4 & 108864 & .&\nn .&\nn .&\nn .&\nn .&\nn 1&\nn 1&\nn .&\nn .&\nn .&\nn .&\nn 1&\nn 1&\nn .&\nn 1&\nn .&\nn .&\nn .&\nn 1&\nn .&\nn .&\nn .&\nn 1&\nn 1 \\[\sss] 
   93 & 4 & 217728 & .&\nn .&\nn .&\nn .&\nn .&\nn .&\nn .&\nn 1&\nn .&\nn .&\nn .&\nn .&\nn 1&\nn .&\nn .&\nn .&\nn .&\nn .&\nn 1&\nn 1&\nn .&\nn 1&\nn .&\nn . \\[\sss] 
   94 & 4 & 217728 & .&\nn .&\nn .&\nn .&\nn .&\nn .&\nn .&\nn 1&\nn .&\nn .&\nn .&\nn 1&\nn 1&\nn .&\nn .&\nn .&\nn 1&\nn .&\nn 1&\nn .&\nn .&\nn .&\nn .&\nn . \\[\sss] 
   95 & 4 & 217728 & .&\nn .&\nn .&\nn .&\nn .&\nn .&\nn .&\nn 1&\nn .&\nn .&\nn .&\nn 1&\nn 1&\nn .&\nn .&\nn .&\nn 1&\nn .&\nn 1&\nn .&\nn .&\nn .&\nn 1&\nn . \\[\sss] 
   96 & 4 & 217728 & .&\nn .&\nn .&\nn .&\nn .&\nn .&\nn .&\nn 1&\nn .&\nn .&\nn 1&\nn .&\nn 1&\nn .&\nn .&\nn .&\nn 1&\nn .&\nn .&\nn .&\nn .&\nn 1&\nn .&\nn . \\[\sss] 
   97 & 4 & 217728 & .&\nn .&\nn .&\nn .&\nn .&\nn .&\nn .&\nn 1&\nn .&\nn .&\nn 1&\nn .&\nn 1&\nn .&\nn .&\nn .&\nn 1&\nn 1&\nn .&\nn .&\nn .&\nn 1&\nn .&\nn . \\[\sss] 
   98 & 4 & 217728 & .&\nn .&\nn .&\nn .&\nn .&\nn .&\nn .&\nn 1&\nn .&\nn .&\nn 1&\nn .&\nn 1&\nn 1&\nn .&\nn .&\nn 1&\nn .&\nn .&\nn .&\nn .&\nn .&\nn 1&\nn . \\[\sss] 
   99 & 4 & 217728 & .&\nn .&\nn .&\nn .&\nn .&\nn .&\nn .&\nn 1&\nn .&\nn .&\nn 1&\nn .&\nn 1&\nn 1&\nn .&\nn .&\nn 1&\nn .&\nn .&\nn .&\nn 1&\nn .&\nn 1&\nn . \\ \midrule 
  100 & 5 &   2268 & .&\nn .&\nn .&\nn .&\nn .&\nn .&\nn .&\nn .&\nn .&\nn .&\nn .&\nn .&\nn .&\nn 1&\nn 1&\nn .&\nn .&\nn 1&\nn 1&\nn .&\nn 1&\nn .&\nn 1&\nn 1 \\[\sss] 
  101 & 5 &   2268 & .&\nn .&\nn .&\nn .&\nn .&\nn .&\nn .&\nn .&\nn .&\nn .&\nn .&\nn 1&\nn .&\nn 1&\nn 1&\nn .&\nn .&\nn .&\nn 1&\nn .&\nn 1&\nn .&\nn 1&\nn 1 \\[\sss] 
  102 & 5 &   4536 & .&\nn .&\nn .&\nn .&\nn .&\nn .&\nn .&\nn .&\nn .&\nn .&\nn .&\nn 1&\nn .&\nn 1&\nn 1&\nn .&\nn .&\nn 1&\nn 1&\nn .&\nn 1&\nn .&\nn .&\nn . \\[\sss] 
  103 & 5 &   9072 & .&\nn .&\nn .&\nn .&\nn .&\nn .&\nn .&\nn .&\nn .&\nn .&\nn .&\nn 1&\nn .&\nn 1&\nn 1&\nn .&\nn .&\nn 1&\nn 1&\nn .&\nn 1&\nn .&\nn .&\nn 1 \\[\sss] 
  104 & 5 &   9072 & .&\nn .&\nn .&\nn .&\nn .&\nn .&\nn .&\nn .&\nn .&\nn .&\nn .&\nn 1&\nn .&\nn 1&\nn 1&\nn .&\nn 1&\nn .&\nn .&\nn .&\nn .&\nn .&\nn .&\nn 1 \\[\sss] 
  105 & 5 &   9072 & .&\nn .&\nn .&\nn .&\nn .&\nn 1&\nn 1&\nn .&\nn .&\nn .&\nn .&\nn .&\nn 1&\nn .&\nn 1&\nn 1&\nn .&\nn 1&\nn 1&\nn .&\nn .&\nn .&\nn .&\nn . \\ \midrule
\end{tabular}
\caption{Large orbits of $3 \times 2 \times 2 \times 2$ tensors, part 1}
\label{table3222part1}
\end{table}  
     
\begin{table} \tiny
\begin{tabular}{rrrcccccccccccccccccccccccc}
$\#$ & rank & size & \multicolumn{24}{c}{canonical form} \\ \midrule  
  106 & 5 &  18144 & .&\nn .&\nn .&\nn .&\nn .&\nn 1&\nn 1&\nn .&\nn .&\nn .&\nn .&\nn .&\nn 1&\nn .&\nn 1&\nn 1&\nn .&\nn 1&\nn 1&\nn 1&\nn .&\nn .&\nn .&\nn . \\[\sss] 
  107 & 5 &  18144 & .&\nn .&\nn .&\nn .&\nn .&\nn 1&\nn 1&\nn .&\nn .&\nn .&\nn .&\nn .&\nn 1&\nn .&\nn 1&\nn 1&\nn .&\nn 1&\nn 1&\nn 1&\nn .&\nn .&\nn 1&\nn . \\[\sss] 
  108 & 5 &  18144 & .&\nn .&\nn .&\nn .&\nn .&\nn 1&\nn 1&\nn .&\nn .&\nn .&\nn .&\nn 1&\nn .&\nn .&\nn 1&\nn .&\nn .&\nn 1&\nn 1&\nn .&\nn 1&\nn .&\nn 1&\nn . \\[\sss] 
  109 & 5 &  18144 & .&\nn .&\nn .&\nn .&\nn .&\nn 1&\nn 1&\nn .&\nn .&\nn .&\nn .&\nn 1&\nn .&\nn .&\nn 1&\nn .&\nn .&\nn 1&\nn 1&\nn .&\nn 1&\nn .&\nn 1&\nn 1 \\[\sss] 
  110 & 5 &  18144 & .&\nn .&\nn .&\nn .&\nn .&\nn 1&\nn 1&\nn .&\nn .&\nn .&\nn .&\nn 1&\nn 1&\nn .&\nn .&\nn .&\nn .&\nn 1&\nn 1&\nn .&\nn .&\nn .&\nn .&\nn . \\[\sss] 
  111 & 5 &  27216 & .&\nn .&\nn .&\nn .&\nn .&\nn .&\nn .&\nn 1&\nn .&\nn .&\nn .&\nn .&\nn .&\nn 1&\nn 1&\nn .&\nn .&\nn 1&\nn 1&\nn .&\nn .&\nn .&\nn .&\nn . \\[\sss] 
  112 & 5 &  27216 & .&\nn .&\nn .&\nn .&\nn .&\nn .&\nn .&\nn 1&\nn .&\nn .&\nn .&\nn .&\nn .&\nn 1&\nn 1&\nn .&\nn .&\nn 1&\nn 1&\nn .&\nn .&\nn .&\nn 1&\nn . \\[\sss] 
  113 & 5 &  27216 & .&\nn .&\nn .&\nn .&\nn .&\nn .&\nn .&\nn 1&\nn .&\nn .&\nn .&\nn .&\nn .&\nn 1&\nn 1&\nn .&\nn .&\nn 1&\nn 1&\nn .&\nn 1&\nn .&\nn .&\nn . \\[\sss] 
  114 & 5 &  27216 & .&\nn .&\nn .&\nn .&\nn .&\nn .&\nn .&\nn 1&\nn .&\nn .&\nn .&\nn .&\nn .&\nn 1&\nn 1&\nn .&\nn .&\nn 1&\nn 1&\nn .&\nn 1&\nn .&\nn 1&\nn . \\[\sss] 
  115 & 5 &  27216 & .&\nn .&\nn .&\nn .&\nn .&\nn .&\nn .&\nn 1&\nn .&\nn .&\nn .&\nn .&\nn 1&\nn .&\nn .&\nn .&\nn 1&\nn .&\nn .&\nn 1&\nn .&\nn 1&\nn 1&\nn . \\[\sss] 
  116 & 5 &  27216 & .&\nn .&\nn .&\nn .&\nn .&\nn .&\nn .&\nn 1&\nn .&\nn .&\nn .&\nn 1&\nn .&\nn 1&\nn 1&\nn .&\nn .&\nn .&\nn 1&\nn .&\nn 1&\nn .&\nn .&\nn . \\[\sss] 
  117 & 5 &  27216 & .&\nn .&\nn .&\nn .&\nn .&\nn .&\nn .&\nn 1&\nn .&\nn .&\nn .&\nn 1&\nn .&\nn 1&\nn 1&\nn .&\nn .&\nn .&\nn 1&\nn .&\nn 1&\nn 1&\nn .&\nn . \\[\sss] 
  118 & 5 &  27216 & .&\nn .&\nn .&\nn .&\nn .&\nn 1&\nn 1&\nn .&\nn .&\nn .&\nn .&\nn .&\nn 1&\nn .&\nn 1&\nn 1&\nn .&\nn 1&\nn 1&\nn .&\nn .&\nn .&\nn .&\nn 1 \\[\sss] 
  119 & 5 &  27216 & .&\nn .&\nn .&\nn .&\nn .&\nn 1&\nn 1&\nn .&\nn .&\nn 1&\nn 1&\nn 1&\nn .&\nn .&\nn .&\nn .&\nn 1&\nn .&\nn .&\nn 1&\nn 1&\nn .&\nn 1&\nn 1 \\[\sss] 
  120 & 5 &  36288 & .&\nn .&\nn .&\nn .&\nn .&\nn .&\nn .&\nn 1&\nn .&\nn .&\nn .&\nn 1&\nn .&\nn 1&\nn 1&\nn .&\nn 1&\nn .&\nn .&\nn .&\nn .&\nn .&\nn .&\nn . \\[\sss] 
  121 & 5 &  36288 & .&\nn .&\nn .&\nn .&\nn .&\nn 1&\nn 1&\nn .&\nn .&\nn .&\nn .&\nn 1&\nn .&\nn .&\nn 1&\nn .&\nn .&\nn 1&\nn 1&\nn .&\nn 1&\nn .&\nn .&\nn . \\[\sss] 
  122 & 5 &  36288 & .&\nn .&\nn .&\nn .&\nn .&\nn 1&\nn 1&\nn .&\nn .&\nn .&\nn .&\nn 1&\nn .&\nn .&\nn 1&\nn .&\nn 1&\nn .&\nn .&\nn .&\nn .&\nn .&\nn .&\nn 1 \\[\sss] 
  123 & 5 &  54432 & .&\nn .&\nn .&\nn .&\nn .&\nn .&\nn .&\nn 1&\nn .&\nn .&\nn 1&\nn .&\nn 1&\nn .&\nn .&\nn .&\nn 1&\nn .&\nn .&\nn 1&\nn 1&\nn 1&\nn .&\nn . \\[\sss] 
  124 & 5 &  54432 & .&\nn .&\nn .&\nn .&\nn .&\nn .&\nn .&\nn 1&\nn .&\nn .&\nn 1&\nn .&\nn 1&\nn .&\nn .&\nn .&\nn 1&\nn .&\nn .&\nn 1&\nn 1&\nn 1&\nn 1&\nn . \\[\sss] 
  125 & 5 &  54432 & .&\nn .&\nn .&\nn .&\nn .&\nn .&\nn .&\nn 1&\nn .&\nn 1&\nn 1&\nn .&\nn 1&\nn .&\nn 1&\nn .&\nn 1&\nn .&\nn 1&\nn .&\nn 1&\nn 1&\nn .&\nn . \\[\sss] 
  126 & 5 &  54432 & .&\nn .&\nn .&\nn .&\nn .&\nn 1&\nn 1&\nn .&\nn .&\nn .&\nn .&\nn 1&\nn .&\nn .&\nn 1&\nn .&\nn .&\nn .&\nn 1&\nn .&\nn 1&\nn .&\nn .&\nn 1 \\[\sss] 
  127 & 5 &  54432 & .&\nn .&\nn .&\nn .&\nn .&\nn 1&\nn 1&\nn .&\nn .&\nn .&\nn .&\nn 1&\nn 1&\nn .&\nn .&\nn .&\nn .&\nn 1&\nn 1&\nn .&\nn 1&\nn .&\nn 1&\nn . \\[\sss] 
  128 & 5 &  54432 & .&\nn .&\nn .&\nn .&\nn .&\nn 1&\nn 1&\nn .&\nn .&\nn .&\nn .&\nn 1&\nn 1&\nn .&\nn .&\nn .&\nn 1&\nn .&\nn .&\nn .&\nn .&\nn .&\nn .&\nn 1 \\[\sss] 
  129 & 5 &  54432 & .&\nn .&\nn .&\nn .&\nn .&\nn 1&\nn 1&\nn .&\nn .&\nn .&\nn .&\nn 1&\nn 1&\nn .&\nn .&\nn .&\nn 1&\nn .&\nn .&\nn .&\nn 1&\nn .&\nn .&\nn 1 \\[\sss] 
  130 & 5 &  54432 & .&\nn .&\nn .&\nn .&\nn .&\nn 1&\nn 1&\nn .&\nn .&\nn .&\nn .&\nn 1&\nn 1&\nn .&\nn 1&\nn .&\nn .&\nn 1&\nn 1&\nn .&\nn .&\nn .&\nn .&\nn . \\[\sss] 
  131 & 5 &  54432 & .&\nn .&\nn .&\nn .&\nn .&\nn 1&\nn 1&\nn .&\nn .&\nn .&\nn .&\nn 1&\nn 1&\nn .&\nn 1&\nn .&\nn .&\nn 1&\nn 1&\nn .&\nn .&\nn .&\nn .&\nn 1 \\[\sss] 
  132 & 5 &  54432 & .&\nn .&\nn .&\nn .&\nn .&\nn 1&\nn 1&\nn .&\nn .&\nn .&\nn .&\nn 1&\nn 1&\nn .&\nn 1&\nn .&\nn .&\nn 1&\nn 1&\nn .&\nn 1&\nn .&\nn 1&\nn 1 \\[\sss] 
  133 & 5 &  54432 & .&\nn .&\nn .&\nn .&\nn .&\nn 1&\nn 1&\nn .&\nn .&\nn .&\nn .&\nn 1&\nn 1&\nn .&\nn 1&\nn .&\nn 1&\nn .&\nn .&\nn .&\nn .&\nn .&\nn 1&\nn 1 \\[\sss] 
  134 & 5 &  54432 & .&\nn .&\nn .&\nn .&\nn .&\nn 1&\nn 1&\nn .&\nn .&\nn .&\nn .&\nn 1&\nn 1&\nn .&\nn 1&\nn .&\nn 1&\nn .&\nn .&\nn .&\nn 1&\nn .&\nn 1&\nn 1 \\[\sss] 
  135 & 5 &  54432 & .&\nn .&\nn .&\nn .&\nn .&\nn 1&\nn 1&\nn .&\nn .&\nn .&\nn 1&\nn .&\nn 1&\nn .&\nn .&\nn 1&\nn .&\nn 1&\nn .&\nn .&\nn 1&\nn .&\nn 1&\nn 1 \\[\sss] 
  136 & 5 &  54432 & .&\nn .&\nn .&\nn .&\nn .&\nn 1&\nn 1&\nn .&\nn .&\nn .&\nn 1&\nn .&\nn 1&\nn .&\nn .&\nn 1&\nn 1&\nn .&\nn .&\nn 1&\nn 1&\nn .&\nn .&\nn 1 \\[\sss] 
  137 & 5 &  54432 & .&\nn .&\nn .&\nn .&\nn .&\nn 1&\nn 1&\nn .&\nn .&\nn 1&\nn 1&\nn 1&\nn .&\nn .&\nn 1&\nn .&\nn 1&\nn .&\nn .&\nn 1&\nn .&\nn .&\nn 1&\nn 1 \\[\sss] 
  138 & 5 &  54432 & .&\nn .&\nn .&\nn .&\nn .&\nn 1&\nn 1&\nn .&\nn .&\nn 1&\nn 1&\nn 1&\nn 1&\nn .&\nn .&\nn .&\nn 1&\nn .&\nn .&\nn 1&\nn 1&\nn .&\nn 1&\nn 1 \\[\sss] 
  139 & 5 &  72576 & .&\nn .&\nn .&\nn 1&\nn .&\nn 1&\nn 1&\nn .&\nn .&\nn .&\nn 1&\nn .&\nn 1&\nn 1&\nn .&\nn .&\nn .&\nn 1&\nn .&\nn .&\nn 1&\nn 1&\nn 1&\nn . \\[\sss] 
  140 & 5 & 108864 & .&\nn .&\nn .&\nn .&\nn .&\nn .&\nn .&\nn 1&\nn .&\nn .&\nn .&\nn .&\nn .&\nn 1&\nn 1&\nn .&\nn 1&\nn .&\nn .&\nn 1&\nn .&\nn .&\nn .&\nn . \\[\sss] 
  141 & 5 & 108864 & .&\nn .&\nn .&\nn .&\nn .&\nn .&\nn .&\nn 1&\nn .&\nn .&\nn .&\nn .&\nn .&\nn 1&\nn 1&\nn .&\nn 1&\nn .&\nn .&\nn 1&\nn .&\nn .&\nn 1&\nn . \\[\sss] 
  142 & 5 & 108864 & .&\nn .&\nn .&\nn .&\nn .&\nn .&\nn .&\nn 1&\nn .&\nn .&\nn .&\nn 1&\nn .&\nn 1&\nn 1&\nn .&\nn 1&\nn .&\nn .&\nn .&\nn .&\nn .&\nn 1&\nn . \\[\sss] 
  143 & 5 & 108864 & .&\nn .&\nn .&\nn .&\nn .&\nn .&\nn .&\nn 1&\nn .&\nn .&\nn .&\nn 1&\nn 1&\nn .&\nn .&\nn .&\nn .&\nn .&\nn 1&\nn .&\nn .&\nn 1&\nn .&\nn . \\[\sss] 
  144 & 5 & 108864 & .&\nn .&\nn .&\nn .&\nn .&\nn .&\nn .&\nn 1&\nn .&\nn .&\nn .&\nn 1&\nn 1&\nn .&\nn .&\nn .&\nn .&\nn .&\nn 1&\nn .&\nn 1&\nn 1&\nn .&\nn . \\[\sss] 
  145 & 5 & 108864 & .&\nn .&\nn .&\nn .&\nn .&\nn .&\nn .&\nn 1&\nn .&\nn .&\nn .&\nn 1&\nn 1&\nn .&\nn .&\nn .&\nn .&\nn 1&\nn 1&\nn .&\nn .&\nn .&\nn .&\nn . \\[\sss] 
  146 & 5 & 108864 & .&\nn .&\nn .&\nn .&\nn .&\nn .&\nn .&\nn 1&\nn .&\nn .&\nn .&\nn 1&\nn 1&\nn .&\nn .&\nn .&\nn .&\nn 1&\nn 1&\nn .&\nn .&\nn .&\nn 1&\nn . \\[\sss] 
  147 & 5 & 108864 & .&\nn .&\nn .&\nn .&\nn .&\nn .&\nn .&\nn 1&\nn .&\nn .&\nn .&\nn 1&\nn 1&\nn .&\nn .&\nn .&\nn .&\nn 1&\nn 1&\nn .&\nn 1&\nn .&\nn .&\nn . \\[\sss] 
  148 & 5 & 108864 & .&\nn .&\nn .&\nn .&\nn .&\nn .&\nn .&\nn 1&\nn .&\nn .&\nn .&\nn 1&\nn 1&\nn .&\nn .&\nn .&\nn .&\nn 1&\nn 1&\nn .&\nn 1&\nn .&\nn 1&\nn . \\[\sss] 
  149 & 5 & 108864 & .&\nn .&\nn .&\nn .&\nn .&\nn .&\nn .&\nn 1&\nn .&\nn .&\nn .&\nn 1&\nn 1&\nn .&\nn .&\nn .&\nn 1&\nn .&\nn .&\nn .&\nn .&\nn 1&\nn 1&\nn . \\[\sss] 
  150 & 5 & 108864 & .&\nn .&\nn .&\nn .&\nn .&\nn .&\nn .&\nn 1&\nn .&\nn .&\nn 1&\nn .&\nn 1&\nn .&\nn .&\nn .&\nn 1&\nn .&\nn .&\nn .&\nn 1&\nn 1&\nn 1&\nn . \\[\sss] 
  151 & 5 & 108864 & .&\nn .&\nn .&\nn .&\nn .&\nn .&\nn .&\nn 1&\nn .&\nn .&\nn 1&\nn .&\nn 1&\nn .&\nn .&\nn .&\nn 1&\nn .&\nn .&\nn 1&\nn .&\nn 1&\nn .&\nn . \\[\sss] 
  152 & 5 & 108864 & .&\nn .&\nn .&\nn .&\nn .&\nn .&\nn .&\nn 1&\nn .&\nn .&\nn 1&\nn .&\nn 1&\nn .&\nn .&\nn .&\nn 1&\nn 1&\nn .&\nn .&\nn .&\nn .&\nn 1&\nn . \\[\sss] 
  153 & 5 & 108864 & .&\nn .&\nn .&\nn .&\nn .&\nn .&\nn .&\nn 1&\nn .&\nn .&\nn 1&\nn .&\nn 1&\nn .&\nn .&\nn .&\nn 1&\nn 1&\nn .&\nn .&\nn 1&\nn 1&\nn 1&\nn . \\[\sss] 
  154 & 5 & 108864 & .&\nn .&\nn .&\nn .&\nn .&\nn .&\nn .&\nn 1&\nn .&\nn 1&\nn 1&\nn .&\nn 1&\nn .&\nn .&\nn .&\nn 1&\nn .&\nn .&\nn .&\nn .&\nn .&\nn 1&\nn . \\[\sss] 
  155 & 5 & 108864 & .&\nn .&\nn .&\nn .&\nn .&\nn .&\nn .&\nn 1&\nn .&\nn 1&\nn 1&\nn .&\nn 1&\nn .&\nn .&\nn .&\nn 1&\nn .&\nn .&\nn .&\nn 1&\nn 1&\nn 1&\nn . \\[\sss] 
  156 & 5 & 108864 & .&\nn .&\nn .&\nn .&\nn .&\nn .&\nn .&\nn 1&\nn .&\nn 1&\nn 1&\nn .&\nn 1&\nn .&\nn 1&\nn .&\nn 1&\nn .&\nn .&\nn .&\nn .&\nn .&\nn 1&\nn . \\[\sss] 
  157 & 5 & 108864 & .&\nn .&\nn .&\nn .&\nn .&\nn .&\nn .&\nn 1&\nn .&\nn 1&\nn 1&\nn .&\nn 1&\nn .&\nn 1&\nn .&\nn 1&\nn .&\nn .&\nn .&\nn .&\nn 1&\nn 1&\nn . \\[\sss] 
  158 & 5 & 108864 & .&\nn .&\nn .&\nn .&\nn .&\nn 1&\nn 1&\nn .&\nn .&\nn .&\nn .&\nn 1&\nn .&\nn .&\nn 1&\nn .&\nn .&\nn .&\nn 1&\nn .&\nn 1&\nn .&\nn 1&\nn . \\[\sss] 
  159 & 5 & 108864 & .&\nn .&\nn .&\nn .&\nn .&\nn 1&\nn 1&\nn .&\nn .&\nn .&\nn .&\nn 1&\nn .&\nn .&\nn 1&\nn .&\nn 1&\nn .&\nn .&\nn .&\nn 1&\nn .&\nn .&\nn 1 \\[\sss] 
  160 & 5 & 108864 & .&\nn .&\nn .&\nn .&\nn .&\nn 1&\nn 1&\nn .&\nn .&\nn .&\nn .&\nn 1&\nn 1&\nn .&\nn .&\nn .&\nn .&\nn .&\nn 1&\nn .&\nn .&\nn .&\nn .&\nn 1 \\[\sss] 
  161 & 5 & 108864 & .&\nn .&\nn .&\nn .&\nn .&\nn 1&\nn 1&\nn .&\nn .&\nn .&\nn .&\nn 1&\nn 1&\nn .&\nn .&\nn .&\nn .&\nn 1&\nn 1&\nn .&\nn .&\nn .&\nn 1&\nn . \\[\sss] 
  162 & 5 & 108864 & .&\nn .&\nn .&\nn .&\nn .&\nn 1&\nn 1&\nn .&\nn .&\nn .&\nn .&\nn 1&\nn 1&\nn .&\nn .&\nn .&\nn .&\nn 1&\nn 1&\nn .&\nn 1&\nn .&\nn .&\nn . \\[\sss] 
  163 & 5 & 108864 & .&\nn .&\nn .&\nn .&\nn .&\nn 1&\nn 1&\nn .&\nn .&\nn .&\nn .&\nn 1&\nn 1&\nn .&\nn .&\nn .&\nn 1&\nn .&\nn .&\nn .&\nn .&\nn .&\nn 1&\nn 1 \\[\sss] 
  164 & 5 & 108864 & .&\nn .&\nn .&\nn .&\nn .&\nn 1&\nn 1&\nn .&\nn .&\nn .&\nn .&\nn 1&\nn 1&\nn .&\nn .&\nn .&\nn 1&\nn .&\nn .&\nn .&\nn 1&\nn .&\nn 1&\nn 1 \\[\sss] 
  165 & 5 & 108864 & .&\nn .&\nn .&\nn .&\nn .&\nn 1&\nn 1&\nn .&\nn .&\nn .&\nn .&\nn 1&\nn 1&\nn .&\nn 1&\nn .&\nn .&\nn 1&\nn 1&\nn .&\nn .&\nn .&\nn 1&\nn . \\[\sss] 
  166 & 5 & 108864 & .&\nn .&\nn .&\nn .&\nn .&\nn 1&\nn 1&\nn .&\nn .&\nn .&\nn .&\nn 1&\nn 1&\nn .&\nn 1&\nn .&\nn .&\nn 1&\nn 1&\nn .&\nn 1&\nn .&\nn .&\nn . \\[\sss] 
  167 & 5 & 108864 & .&\nn .&\nn .&\nn .&\nn .&\nn 1&\nn 1&\nn .&\nn .&\nn .&\nn .&\nn 1&\nn 1&\nn .&\nn 1&\nn .&\nn 1&\nn .&\nn 1&\nn .&\nn .&\nn .&\nn .&\nn 1 \\[\sss] 
  168 & 5 & 108864 & .&\nn .&\nn .&\nn .&\nn .&\nn 1&\nn 1&\nn .&\nn .&\nn .&\nn .&\nn 1&\nn 1&\nn .&\nn 1&\nn .&\nn 1&\nn .&\nn 1&\nn .&\nn 1&\nn .&\nn 1&\nn 1 \\[\sss] 
  169 & 5 & 108864 & .&\nn .&\nn .&\nn .&\nn .&\nn 1&\nn 1&\nn .&\nn .&\nn .&\nn 1&\nn .&\nn .&\nn .&\nn .&\nn 1&\nn .&\nn 1&\nn .&\nn .&\nn 1&\nn .&\nn .&\nn . \\[\sss] 
  170 & 5 & 108864 & .&\nn .&\nn .&\nn .&\nn .&\nn 1&\nn 1&\nn .&\nn .&\nn .&\nn 1&\nn .&\nn .&\nn .&\nn .&\nn 1&\nn .&\nn 1&\nn .&\nn .&\nn 1&\nn .&\nn 1&\nn . \\[\sss] 
  171 & 5 & 108864 & .&\nn .&\nn .&\nn .&\nn .&\nn 1&\nn 1&\nn .&\nn .&\nn .&\nn 1&\nn .&\nn 1&\nn .&\nn .&\nn 1&\nn .&\nn 1&\nn .&\nn 1&\nn 1&\nn .&\nn .&\nn . \\[\sss] 
  172 & 5 & 108864 & .&\nn .&\nn .&\nn .&\nn .&\nn 1&\nn 1&\nn .&\nn .&\nn .&\nn 1&\nn .&\nn 1&\nn .&\nn .&\nn 1&\nn .&\nn 1&\nn .&\nn 1&\nn 1&\nn .&\nn 1&\nn 1 \\[\sss] 
  173 & 5 & 108864 & .&\nn .&\nn .&\nn .&\nn .&\nn 1&\nn 1&\nn .&\nn .&\nn .&\nn 1&\nn .&\nn 1&\nn .&\nn .&\nn 1&\nn 1&\nn .&\nn .&\nn 1&\nn .&\nn .&\nn 1&\nn . \\[\sss] 
  174 & 5 & 108864 & .&\nn .&\nn .&\nn .&\nn .&\nn 1&\nn 1&\nn .&\nn .&\nn 1&\nn 1&\nn .&\nn .&\nn .&\nn .&\nn 1&\nn 1&\nn .&\nn 1&\nn 1&\nn .&\nn .&\nn 1&\nn . \\[\sss] 
  175 & 5 & 108864 & .&\nn .&\nn .&\nn .&\nn .&\nn 1&\nn 1&\nn .&\nn .&\nn 1&\nn 1&\nn .&\nn .&\nn .&\nn 1&\nn .&\nn 1&\nn .&\nn 1&\nn 1&\nn .&\nn .&\nn .&\nn 1 \\[\sss] 
  176 & 5 & 108864 & .&\nn .&\nn .&\nn .&\nn .&\nn 1&\nn 1&\nn .&\nn .&\nn 1&\nn 1&\nn 1&\nn .&\nn .&\nn 1&\nn .&\nn 1&\nn .&\nn .&\nn 1&\nn .&\nn .&\nn .&\nn 1 \\[\sss] 
  177 & 5 & 217728 & .&\nn .&\nn .&\nn .&\nn .&\nn .&\nn .&\nn 1&\nn .&\nn .&\nn .&\nn 1&\nn 1&\nn .&\nn .&\nn .&\nn 1&\nn .&\nn .&\nn .&\nn .&\nn .&\nn 1&\nn . \\[\sss] 
  178 & 5 & 217728 & .&\nn .&\nn .&\nn .&\nn .&\nn .&\nn .&\nn 1&\nn .&\nn .&\nn .&\nn 1&\nn 1&\nn .&\nn .&\nn .&\nn 1&\nn .&\nn 1&\nn .&\nn .&\nn 1&\nn .&\nn . \\[\sss] 
  179 & 5 & 217728 & .&\nn .&\nn .&\nn .&\nn .&\nn .&\nn .&\nn 1&\nn .&\nn .&\nn .&\nn 1&\nn 1&\nn .&\nn .&\nn .&\nn 1&\nn .&\nn 1&\nn .&\nn .&\nn 1&\nn 1&\nn . \\[\sss] 
  180 & 5 & 217728 & .&\nn .&\nn .&\nn .&\nn .&\nn .&\nn .&\nn 1&\nn .&\nn .&\nn .&\nn 1&\nn 1&\nn .&\nn .&\nn .&\nn 1&\nn 1&\nn 1&\nn .&\nn .&\nn .&\nn 1&\nn . \\[\sss] 
  181 & 5 & 217728 & .&\nn .&\nn .&\nn .&\nn .&\nn .&\nn .&\nn 1&\nn .&\nn .&\nn 1&\nn .&\nn 1&\nn 1&\nn .&\nn .&\nn 1&\nn .&\nn .&\nn 1&\nn .&\nn .&\nn .&\nn . \\[\sss] 
  182 & 5 & 217728 & .&\nn .&\nn .&\nn .&\nn .&\nn .&\nn .&\nn 1&\nn .&\nn .&\nn 1&\nn .&\nn 1&\nn 1&\nn .&\nn .&\nn 1&\nn .&\nn .&\nn 1&\nn .&\nn .&\nn 1&\nn . \\[\sss] 
  183 & 5 & 217728 & .&\nn .&\nn .&\nn .&\nn .&\nn .&\nn .&\nn 1&\nn .&\nn .&\nn 1&\nn .&\nn 1&\nn 1&\nn .&\nn .&\nn 1&\nn .&\nn .&\nn 1&\nn 1&\nn .&\nn .&\nn . \\[\sss] 
  184 & 5 & 217728 & .&\nn .&\nn .&\nn .&\nn .&\nn .&\nn .&\nn 1&\nn .&\nn .&\nn 1&\nn .&\nn 1&\nn 1&\nn .&\nn .&\nn 1&\nn .&\nn .&\nn 1&\nn 1&\nn .&\nn 1&\nn . \\[\sss] 
  185 & 5 & 217728 & .&\nn .&\nn .&\nn .&\nn .&\nn .&\nn .&\nn 1&\nn .&\nn 1&\nn 1&\nn .&\nn 1&\nn .&\nn .&\nn .&\nn 1&\nn .&\nn .&\nn .&\nn 1&\nn .&\nn 1&\nn . \\[\sss] 
  186 & 5 & 217728 & .&\nn .&\nn .&\nn .&\nn .&\nn 1&\nn 1&\nn .&\nn .&\nn .&\nn .&\nn 1&\nn 1&\nn .&\nn .&\nn .&\nn .&\nn .&\nn 1&\nn .&\nn .&\nn .&\nn 1&\nn 1 \\[\sss] 
  187 & 5 & 217728 & .&\nn .&\nn .&\nn .&\nn .&\nn 1&\nn 1&\nn .&\nn .&\nn .&\nn .&\nn 1&\nn 1&\nn .&\nn .&\nn .&\nn 1&\nn .&\nn 1&\nn .&\nn .&\nn .&\nn .&\nn 1 \\[\sss] 
  188 & 5 & 217728 & .&\nn .&\nn .&\nn .&\nn .&\nn 1&\nn 1&\nn .&\nn .&\nn .&\nn .&\nn 1&\nn 1&\nn .&\nn .&\nn .&\nn 1&\nn .&\nn 1&\nn .&\nn .&\nn .&\nn 1&\nn . \\[\sss] 
  189 & 5 & 217728 & .&\nn .&\nn .&\nn .&\nn .&\nn 1&\nn 1&\nn .&\nn .&\nn .&\nn .&\nn 1&\nn 1&\nn .&\nn .&\nn .&\nn 1&\nn .&\nn 1&\nn .&\nn .&\nn .&\nn 1&\nn 1 \\[\sss] 
  190 & 5 & 217728 & .&\nn .&\nn .&\nn .&\nn .&\nn 1&\nn 1&\nn .&\nn .&\nn .&\nn .&\nn 1&\nn 1&\nn .&\nn 1&\nn .&\nn 1&\nn .&\nn 1&\nn .&\nn 1&\nn .&\nn .&\nn . \\[\sss] 
  191 & 5 & 217728 & .&\nn .&\nn .&\nn .&\nn .&\nn 1&\nn 1&\nn .&\nn .&\nn .&\nn 1&\nn .&\nn .&\nn .&\nn .&\nn 1&\nn .&\nn 1&\nn .&\nn .&\nn 1&\nn .&\nn .&\nn 1 \\[\sss] 
  192 & 5 & 217728 & .&\nn .&\nn .&\nn .&\nn .&\nn 1&\nn 1&\nn .&\nn .&\nn .&\nn 1&\nn .&\nn .&\nn .&\nn .&\nn 1&\nn .&\nn 1&\nn .&\nn .&\nn 1&\nn .&\nn 1&\nn 1 \\[\sss] 
  193 & 5 & 217728 & .&\nn .&\nn .&\nn .&\nn .&\nn 1&\nn 1&\nn .&\nn .&\nn .&\nn 1&\nn .&\nn .&\nn .&\nn .&\nn 1&\nn .&\nn 1&\nn .&\nn 1&\nn 1&\nn .&\nn .&\nn . \\[\sss] 
  194 & 5 & 217728 & .&\nn .&\nn .&\nn .&\nn .&\nn 1&\nn 1&\nn .&\nn .&\nn .&\nn 1&\nn .&\nn .&\nn .&\nn .&\nn 1&\nn .&\nn 1&\nn .&\nn 1&\nn 1&\nn .&\nn 1&\nn . \\[\sss] 
  195 & 5 & 217728 & .&\nn .&\nn .&\nn .&\nn .&\nn 1&\nn 1&\nn .&\nn .&\nn .&\nn 1&\nn .&\nn .&\nn .&\nn .&\nn 1&\nn .&\nn 1&\nn .&\nn 1&\nn 1&\nn .&\nn 1&\nn 1 \\[\sss] 
  196 & 5 & 217728 & .&\nn .&\nn .&\nn .&\nn .&\nn 1&\nn 1&\nn .&\nn .&\nn .&\nn 1&\nn .&\nn 1&\nn .&\nn .&\nn 1&\nn 1&\nn .&\nn .&\nn 1&\nn .&\nn .&\nn .&\nn 1 \\[\sss] 
  197 & 5 & 217728 & .&\nn .&\nn .&\nn .&\nn .&\nn 1&\nn 1&\nn .&\nn .&\nn .&\nn 1&\nn .&\nn 1&\nn .&\nn .&\nn 1&\nn 1&\nn .&\nn .&\nn 1&\nn .&\nn .&\nn 1&\nn 1 \\[\sss] 
  198 & 5 & 217728 & .&\nn .&\nn .&\nn 1&\nn .&\nn 1&\nn 1&\nn .&\nn .&\nn .&\nn 1&\nn .&\nn 1&\nn 1&\nn .&\nn .&\nn 1&\nn .&\nn .&\nn .&\nn .&\nn 1&\nn .&\nn 1 \\[\sss] 
  199 & 5 & 217728 & .&\nn .&\nn .&\nn 1&\nn .&\nn 1&\nn 1&\nn .&\nn .&\nn .&\nn 1&\nn .&\nn 1&\nn 1&\nn .&\nn .&\nn 1&\nn .&\nn .&\nn .&\nn 1&\nn .&\nn 1&\nn 1 \\ \midrule
  200 & 6 &    168 & .&\nn .&\nn .&\nn .&\nn .&\nn .&\nn .&\nn .&\nn .&\nn 1&\nn 1&\nn .&\nn 1&\nn .&\nn 1&\nn 1&\nn 1&\nn .&\nn 1&\nn 1&\nn 1&\nn 1&\nn .&\nn 1 \\[\sss] 
  201 & 6 &   6048 & .&\nn .&\nn .&\nn .&\nn .&\nn 1&\nn 1&\nn .&\nn .&\nn 1&\nn 1&\nn .&\nn 1&\nn .&\nn 1&\nn 1&\nn 1&\nn .&\nn 1&\nn 1&\nn .&\nn .&\nn .&\nn . \\[\sss] 
  202 & 6 &  12096 & .&\nn .&\nn .&\nn 1&\nn .&\nn 1&\nn 1&\nn .&\nn .&\nn 1&\nn 1&\nn .&\nn 1&\nn .&\nn .&\nn 1&\nn 1&\nn .&\nn .&\nn .&\nn .&\nn .&\nn .&\nn 1 \\[\sss] 
  203 & 6 &  18144 & .&\nn .&\nn .&\nn .&\nn .&\nn .&\nn .&\nn 1&\nn .&\nn .&\nn .&\nn 1&\nn .&\nn 1&\nn 1&\nn .&\nn .&\nn 1&\nn 1&\nn .&\nn 1&\nn .&\nn .&\nn . \\[\sss] 
  204 & 6 &  18144 & .&\nn .&\nn .&\nn .&\nn .&\nn .&\nn .&\nn 1&\nn .&\nn 1&\nn 1&\nn .&\nn 1&\nn .&\nn 1&\nn .&\nn 1&\nn .&\nn 1&\nn 1&\nn 1&\nn 1&\nn .&\nn . \\[\sss] 
  205 & 6 &  18144 & .&\nn .&\nn .&\nn .&\nn .&\nn 1&\nn 1&\nn .&\nn .&\nn 1&\nn 1&\nn 1&\nn 1&\nn .&\nn .&\nn .&\nn 1&\nn .&\nn .&\nn 1&\nn .&\nn .&\nn .&\nn 1 \\[\sss] 
  206 & 6 &  54432 & .&\nn .&\nn .&\nn .&\nn .&\nn 1&\nn 1&\nn .&\nn .&\nn .&\nn .&\nn 1&\nn 1&\nn .&\nn .&\nn .&\nn .&\nn 1&\nn 1&\nn .&\nn .&\nn .&\nn .&\nn 1 \\[\sss] 
  207 & 6 &  54432 & .&\nn .&\nn .&\nn .&\nn .&\nn 1&\nn 1&\nn .&\nn .&\nn .&\nn .&\nn 1&\nn 1&\nn .&\nn .&\nn .&\nn .&\nn 1&\nn 1&\nn .&\nn 1&\nn .&\nn 1&\nn 1 \\[\sss] 
  208 & 6 &  54432 & .&\nn .&\nn .&\nn .&\nn .&\nn 1&\nn 1&\nn .&\nn .&\nn 1&\nn 1&\nn .&\nn 1&\nn .&\nn 1&\nn 1&\nn 1&\nn .&\nn 1&\nn 1&\nn .&\nn .&\nn .&\nn 1 \\[\sss] 
  209 & 6 &  54432 & .&\nn .&\nn .&\nn 1&\nn .&\nn 1&\nn 1&\nn .&\nn .&\nn .&\nn 1&\nn .&\nn 1&\nn .&\nn 1&\nn 1&\nn .&\nn 1&\nn .&\nn .&\nn 1&\nn 1&\nn 1&\nn . \\[\sss] 
  210 & 6 & 108864 & .&\nn .&\nn .&\nn .&\nn .&\nn 1&\nn 1&\nn .&\nn .&\nn .&\nn 1&\nn .&\nn 1&\nn .&\nn .&\nn 1&\nn 1&\nn .&\nn .&\nn 1&\nn 1&\nn .&\nn 1&\nn 1 \\[\sss] 
  211 & 6 & 108864 & .&\nn .&\nn .&\nn 1&\nn .&\nn 1&\nn 1&\nn .&\nn .&\nn 1&\nn 1&\nn .&\nn 1&\nn .&\nn .&\nn .&\nn 1&\nn .&\nn .&\nn .&\nn .&\nn .&\nn 1&\nn 1 \\[\sss] 
  212 & 6 & 217728 & .&\nn .&\nn .&\nn 1&\nn .&\nn 1&\nn 1&\nn .&\nn .&\nn .&\nn 1&\nn .&\nn 1&\nn 1&\nn .&\nn .&\nn 1&\nn .&\nn .&\nn .&\nn .&\nn 1&\nn 1&\nn 1 \\ \midrule
\end{tabular}
\caption{Large orbits of $3 \times 2 \times 2 \times 2$ tensors, part 2}
\label{table3222part2}
\end{table}  
   
%%%%%%%%%%%%%%%%%%%%%%%%%%%%%%%%%%%%%%%%%%%%%%%%%%%%%%%%%%%%%%%%%%%%%%%%%%%%%%

\section*{Acknowledgement}

The first author was partially supported by a Discovery Grant from NSERC.

%%%%%%%%%%%%%%%%%%%%%%%%%%%%%%%%%%%%%%%%%%%%%%%%%%%%%%%%%%%%%%%%%%%%%%%%%%%%%%

\end{document}